\documentclass[11pt]{article} %
\makeatletter

\newif\ifuserevtex
\newif\ifusetwocolfigs

\@ifclassloaded{revtex4-2}{%
  \userevtextrue
  \@ifclasswith{revtex4-2}{reprint}{%
    \usetwocolfigstrue
  }{%
    \usetwocolfigsfalse
  }%
}{%
  \userevtexfalse
  \usetwocolfigsfalse %
}

\makeatother

\newif\ifshowlinenumbers
\usepackage[T1]{fontenc}         %
\usepackage[utf8]{inputenc}      %
\usepackage[english]{babel}      %

\usepackage{amssymb}             %
\usepackage{mathtools}           %
\usepackage{dcolumn}             %

\usepackage{graphicx}            %
\graphicspath{{./}}

\usepackage{color,xcolor}        %
\usepackage{soul}                %
\usepackage{tcolorbox}           %
\usepackage{cancel}              %
\usepackage{comment}             %
\usepackage{moreverb}            %
\usepackage[normalem]{ulem}      %

\ifuserevtex
    \newcommand{\mybibliostyle}{apsrev4-2}
\else
    \usepackage{amsthm}
    
    \usepackage[a4paper,margin=2.3cm]{geometry} %
    \usepackage{appendix} %
    
    \usepackage[numbers]{natbib} %
    \newcommand{\mybibliostyle}{unsrtnat} %

    \usepackage{authblk} 
    \usepackage{setspace} %
\fi

\usepackage{pifont}

\usepackage{etoolbox}  %
\usepackage{environ}   %
\usepackage{xparse}
\usepackage{scalerel}

\usepackage[
  colorlinks=true,
  linkcolor=blue,  %
  citecolor=blue,  %
  urlcolor=blue    %
]{hyperref}

\usepackage[capitalise]{cleveref}  %

\usepackage{physics}  %

\usepackage{tikz}
\usetikzlibrary{shapes.geometric, arrows.meta, positioning, calc, backgrounds}
\usepackage{algorithm}
\usepackage[noend]{algpseudocode}
\usepackage{listings}
\usepackage{xcolor}

\definecolor{Pcomment}{RGB}{160, 160, 160}   
\definecolor{Pkeyword}{RGB}{0, 136, 0}       
\definecolor{DarkerPastieBlue}{HTML}{007ACC} 
\definecolor{Pstring}{RGB}{186, 33, 33}      
\definecolor{Pbg}{RGB}{255, 255, 255}        

\makeatletter
\def\BState{\State\hskip-\ALG@thistlm}
\makeatother

\newif\ifshowallemails
\showallemailsfalse  %

\makeatletter

\def\affillist{}
\newcommand{\myaffil}[2]{%
  \expandafter\def\csname affil@#1\endcsname{#2}%
  \ifuserevtex\else
    \listgadd{\affillist}{#1}%
  \fi
}

\newcommand{\buildaffiliations}{%
  \ifuserevtex\else
    \forlistloop{\emitoneaffil}{\affillist}
  \fi
}
\newcommand{\emitoneaffil}[1]{%
  \expandafter\affil\expandafter[#1]{\csname affil@#1\endcsname}%
}

\newcommand{\orcidlink}[1]{%
  \href{https://orcid.org/#1}{\textsuperscript{\textcolor{green!50!black}{\textbullet}}}%
}

\NewDocumentCommand{\myauthor}{O{} O{} O{} m m}{%
  \ifuserevtex
    \author{#4%
      \ifx&#2&\else\ \orcidlink{#2}\fi
    }%
    \ifx&#1&\else
      \ifshowallemails
        \email{#1}%
      \else
        \ifx&#3&\else  %
          \email{#1}%
        \fi
      \fi
    \fi
    \@for\aid:=#5\do{%
      \expandafter\affiliation\expandafter{\csname affil@\aid\endcsname}%
    }%
  \else
    \author[#5]{#4%
      \ifx&#2&\else\ \orcidlink{#2}\fi
      \ifx&#1&\else
        \ifx&#3&%
          \ifshowallemails
            \thanks{\href{mailto:#1}{#1}}%
          \fi
        \else
          \thanks{Corresponding author: \href{mailto:#1}{#1}}%
        \fi
      \fi
    }%
  \fi
}

\newcommand{\storedabstract}{}  %
\newcommand{\storedkeywords}{}  %

\newcommand{\myabstract}[1]{%
  \long\def\storedabstract{#1}%
}

\newcommand{\mykeywords}[1]{%
  \def\storedkeywords{#1}%
}

\newcommand{\maketitleandabstract}{%
  \ifuserevtex
    \begin{abstract}
    \storedabstract
    \end{abstract}
    \keywords{\storedkeywords}
    \maketitle
  \else
    \maketitle
    \section*{Abstract}
    \storedabstract
    \ifx\storedkeywords\@empty
    \else
      \vspace{1em}
      \par\noindent \textbf{Keywords:} \storedkeywords
    \fi
  \fi
}

\newcommand{\storedhighlights}{}  %

\newcommand{\myhighlights}[1]{%
  \ifuserevtex
  \else
    \gdef\storedhighlights{#1}%
  \fi
}

\newcommand{\inserthighlights}{%
  \ifuserevtex
  \else
    \par\noindent\textbf{Highlights:}
    \begin{itemize}
      \storedhighlights
    \end{itemize}
    \newpage
  \fi
}

\newcommand{\storedacknowledgements}{}  %

\newcommand{\myacknowledgements}[1]{%
  \gdef\storedacknowledgements{#1}%
}

\newcommand{\insertacknowledgements}{%
  \ifx\storedacknowledgements\@empty
  \else
    \ifuserevtex
      \begin{acknowledgments}
      \storedacknowledgements
      \end{acknowledgments}
    \else
      \section*{Acknowledgements}
      \storedacknowledgements
    \fi
  \fi
}

\makeatother

\ifshowlinenumbers
  \usepackage{lineno}
  \AtBeginDocument{%
    \ifuserevtex
      \ifusetwocolfigs
      \else
        \linenumbers
      \fi
    \else
      \linenumbers
    \fi
  }
\fi

\makeatletter

\def\onecolfig{\@ifnextchar[{\onecolfig@opt}{\onecolfig@opt[]}}
\def\endonecolfig{\end{figure}}
\def\onecolfig@opt[#1]{\begin{figure}[#1]}

\def\twocolfig{\@ifnextchar[{\twocolfig@opt}{\twocolfig@opt[]}}
\def\endtwocolfig{\ifusetwocolfigs \end{figure*} \else \end{figure} \fi}
\def\twocolfig@opt[#1]{%
  \ifusetwocolfigs
    \begin{figure*}[#1]%
  \else
    \begin{figure}[#1]%
  \fi
}

\makeatother

\usepackage{xr}

\makeatletter
\newcommand*{\addFileDependency}[1]{%
  \typeout{(#1)}
  \@addtofilelist{#1}
  \IfFileExists{#1}{}{\typeout{No file #1.}}
}
\makeatother

\newcommand{\papertitle}{ A versatile FEM framework with native GPU scalability via globally-applied AD}

\newcommand{\ie}{{\textit{i.e.}}}
\newcommand{\eg}{{\textit{e.g.}}}

\newcommand{\tatva}{\texttt{tatva}~} 
\begin{document}

\title{\papertitle}
\myaffil{1}{Institute for Building Materials, ETH Zurich, Switzerland}

\myauthor[mpundir@ethz.ch][0000-0000-0000-0000][yes]{Mohit Pundir}{1}

\myauthor[florez@ethz.ch][orcid]{Flavio Lorez}{1}

\myauthor[dkammer@ethz.ch][0000-0003-3782-9368][]{David S. Kammer}{1}

\buildaffiliations \date{\today}

\myabstract{
    Energy-based finite-element formulations provide a unified framework for describing complex physical systems in computational mechanics. In these energy-based methods, the governing equations can be obtained directly by considering the derivatives of a single global energy functional. While Automatic Differentiation (AD) can be used to automate the generation of these derivatives, current frameworks face a clear trade-off based primarily on the scale upon which the AD method is applied. Globally applied AD offers high expressivity but cannot currently be scaled to large problems. Locally applied AD scales well through traditional assembly methods, but the variety of physics and couplings that the framework can easily represent is more limited than the global approach. Here, we introduce an energy-centric framework (\href{https://github.com/smec-ethz/tatva}{tatva}) that defines the physics of a problem as a single global functional and applies AD globally to generate residual and tangent operators. By leveraging Jacobian-vector products for matrix-free solvers and coloring-based sparse differentiation for materializing sparse tangent stiffness matrices when needed, our flexible design scales linearly with the problem size on GPUs. We demonstrate that our framework can handle large problems (with millions of degrees of freedom) without memory exhaustion. Additionally, it offers a unified, fully differentiable methodology that can address a wide range of problems, including multi-point constraints, mixed-dimensional coupling, and the incorporation of neural networks, while maintaining high performance and scalability on modern GPU architectures.
}

\mykeywords{Finite-Element Analysis, Automatic Differentiation, High-Performance Computing, Graph Coloring, GPUs}

\myhighlights{
    \item Introduces a monolithic FEM framework based on global automatic differentiation
    \item Allows ideal linear scalability $\mathcal{O}(N)$ on GPUs 
    \item Replaces rigid element-wise assembly by physics-agnostic global FEM assembly
    \item Enables both matrix-free JVPs and sparse differentiation in FEM
    \item Preserves expressive energy-based formulations while scaling to large problems
}

\maketitleandabstract
\newline
\inserthighlights

\section{Introduction} \label{sec:introduction}

Many variational finite-element methods used in computational mechanics are built on an energy-based formulation~\cite{zienkiewicz_finite_2000, bathe_finite_1996, boffi_mixed_2013}. Based on the principle of minimum potential energy, the governing equations for the equilibrium state are derived as the first- and second-order derivatives of a single scalar energy functional. This approach provides a compact and physically consistent description of complex systems. Representing a problem via an energy functional enables the seamless integration of nonlinear material laws, multiphysics couplings, and constraints within a unified finite-element framework. 

The computation of the required derivatives of the energy functional is often non-trivial for complex problems. Automatic differentiation (AD)~\cite{margossian_review_2019, korelc_automation_2016} offers a systematic way to compute these derivatives. By algorithmically applying the chain rule to complex expressions, including those with conditional logic and iterative loops, AD automatically computes derivatives up to machine precision~\cite{baydin_automatic_2017, rothe_automatic_2015}. Consequently, by delegating the evaluation of derivatives to the computational framework, AD reduces implementation effort and enables rapid exploration of novel physical models. This has led to the widespread adoption of AD in scientific computing~\cite{schoenholz_jax_2020, carrer_learning_2024, toshev_jax-sph_2024, xue_jax-fem_2023, bleyer_numerical_2024, latyshev_expressing_2025}. In principle, AD allows finite-element formulations to be expressed directly at the level of the energy functional.

\begin{figure}
    \centering
    \includegraphics[width=\linewidth]{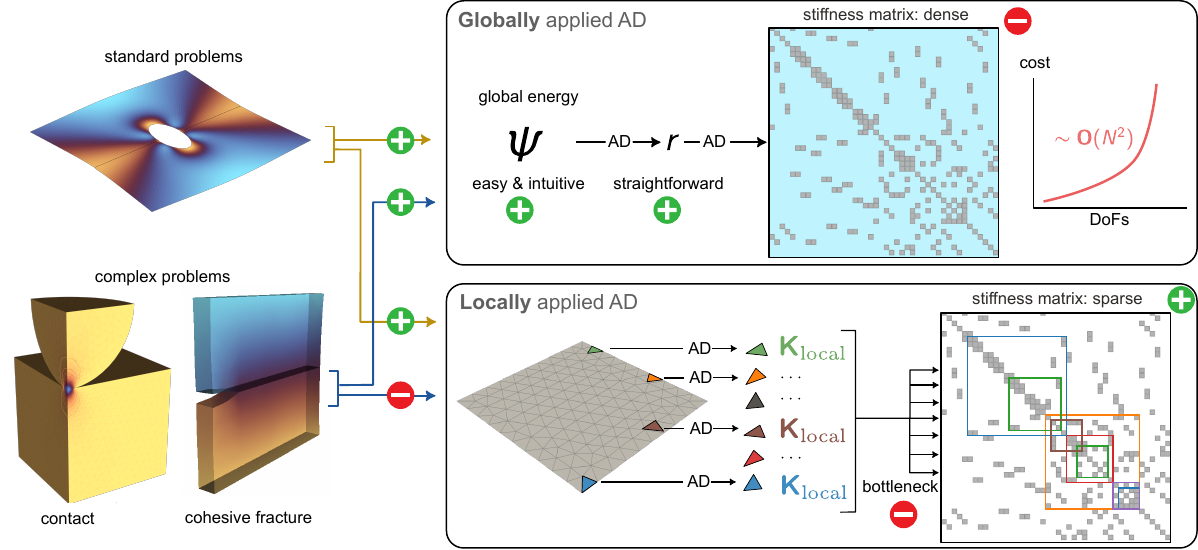}
    \caption{Schematic of global and local approaches in AD-based finite-element frameworks. The global approach differentiates the total energy to obtain residual and tangent operators, enabling complex problems such as contact between deformable bodies and interface fracture, but with quadratic memory and computational costs due to the construction of a dense tangent stiffness matrix. The local approach applies AD at the element level to obtain local tangent operators, then uses standard scatter-add assembly to form global operators, reducing memory cost but making complex problems harder to model and less expressive.}
    \label{fig:scatter-add-workflow}
\end{figure}

In practice, AD has been applied in two distinct ways in finite-element methods. One approach applies AD to a global discrete energy functional~\cite{vigliotti_automatic_2021}. In this setting, the residual vector and tangent operator are obtained by differentiating a single scalar-valued energy with respect to all degrees of freedom, see \Cref{fig:scatter-add-workflow}. This approach provides a direct correspondence between the variational formulation and its discrete realization. It allows for the natural expression of global constraints, nonlocal interactions, mixed-dimensional couplings, and unconventional discretizations. However, direct differentiation of a global energy produces dense derivative operators leading to a memory footprint and computational costs that scale quadratically with the number of degrees of freedom, \ie{} $\mathcal{O}(N^2)$. This overhead typically limits the applicability of the global approach to small-scale problems or specific research applications. 

A second approach restricts AD to local quantities, typically at the element or quadrature-point level~\cite{xue_jax-fem_2023, wu_framework_2023}. Here, AD is used to compute local residuals and tangent contributions, which are then assembled into global sparse operators using standard finite-element procedures, see \Cref{fig:scatter-add-workflow}. This approach addresses memory limitations by preserving the sparsity of the tangent stiffness matrix, thereby enabling scalability to large problem sizes. Consequently, the application of AD to local quantities has become the prevailing approach in practical finite-element implementations that rely on AD~\cite{ichihara_naruki-ichiharafeax_2026}.

However, restricting differentiation to local quantities imposes significant limitations on the simplicity and versatility of the method. By shifting AD to the element level, the global problem structure is no longer directly represented in the differentiation process. Instead, the global system must be recovered through element-wise assembly. As a result, this approach inherits the element-centric abstractions of classical finite-element libraries~\cite{logg_ffc_2012, baratta_dolfinx_2023, anderson_mfem_2021, kaczmarczyk_mofem_2020, richart_akantu_2024}. This makes it difficult to use this framework to analyze systems involving local and non-local constraints (\eg{}, contact or multi-point constraints), mixed-dimensional couplings (\eg{}, cohesive-fracture or beam-solid interaction) or multiphysics systems with disparate degrees of freedom (\eg{}, thermo-mechanical coupling). While AD simplifies the math within an element, it does not fundamentally change how the global problem is formulated or executed. 

The element-level approach also retains the traditional compute–assemble–solve workflow, see \Cref{fig:scatter-add-workflow}. Local contributions are accumulated into global sparse operators via ``scatter-add'' operations that rely on unstructured memory access patterns. While effective on CPUs, this paradigm is not well suited to modern accelerator hardware such as GPUs or TPUs. The irregular memory access required by scatter-add assembly conflicts with the stream-oriented execution model of GPUs. This leads to memory-bandwidth bottlenecks that limit performance and hence scalability. As a result, current research often focuses on optimizing the assembly process~\cite{andrej_scalable_2025, anderson_mfem_2021}; such approaches overlook more substantial improvements that could be obtained by reconsidering the representation of the global variational problem itself.

In summary, existing approaches illustrate a clear trade-off. Global energy-based differentiation offers high expressive power and a transparent connection to the variational formulation, but lacks the computational practicality needed for large-scale simulations. Conversely, element-level differentiation achieves scalability and performance but constrains the class of physical problems that can be expressed naturally. This leaves an unresolved challenge in AD-based finite-element frameworks: achieving both high-level expressiveness to tackle diverse problems and computational efficiency within a single framework.

Here, we revisit global energy-centric formulations in the context of modern AD and high-performance computing tools. We aim to retain the flexibility of global formulations while enabling the efficient evaluation of residuals and tangent operators for large-scale problems. To achieve this, our framework utilizes graph-coloring-based sparse differentiation and Jacobian-vector products (JVPs) to bypass the memory bottlenecks of dense differentiation. We present this approach through \tatva\cite{tatva_library_2026}, a monolithic, energy-centric library. We demonstrate that this framework is both broadly applicable to diverse physical problems and highly scalable on GPU architectures.
 
\section{Proposed framework} \label{sec:proposed_framework}

We propose a framework that adopts an energy-centric philosophy, in which the behavior of a physical system is defined by a scalar potential energy. Following a variational approach~\cite{zienkiewicz_finite_2000}, residuals and tangent operators are obtained as the first and second derivatives of the total energy $\Psi$,
\begin{equation}\label{eq:fundamental-equation}
	\boldsymbol{r}(\boldsymbol{u}) := \nabla \Psi(\boldsymbol{u}), \quad \mathbf{K}(\boldsymbol{u}) := \nabla^2 \Psi(\boldsymbol{u}).
\end{equation}
where $\boldsymbol{u}$ denotes a physical field such as displacement, temperature, or velocity. While this mathematical foundation is standard in variational finite-element formulations (see \eg{}~\cite{vigliotti_automatic_2021}), the proposed framework adopts a different computational architecture. In contrast to existing AD-based approaches~\cite{xue_jax-fem_2023, ichihara_naruki-ichiharafeax_2026, wu_framework_2023}, which embed differentiation within element-level assembly loops, a monolithic structure is introduced here that decouples the description of the physical problem from the numerical backend. %

To realize the monolithic architecture, we discretize the continuous functional $\Psi(\boldsymbol{u})$ as a sum over the finite-element mesh. We evaluate the global potential by aggregating local element contributions computed through numerical integration at quadrature points,
\begin{equation}\label{eq:discretized-energy}
	\Psi_h(\boldsymbol{u}_h) = \sum_{e=1}^{N_{\mathrm{el}}} \sum_{q=1}^{N_q} \psi\!\left(\nabla \boldsymbol{u}_h^{e} (\xi_q)\right)\,  w_{q}\, \det\!\left(\mathbf{J}^{e}(\xi_q) \right)
\end{equation}
where $\psi$ is the local energy density, $w_{q}$ and $\xi_q$ denote the quadrature weights and positions, and $\mathbf{J}^{e}$ is the Jacobian of the mapping from the reference to the physical element.

\begin{figure}[t]
	\centering
	\includegraphics[width=\linewidth]{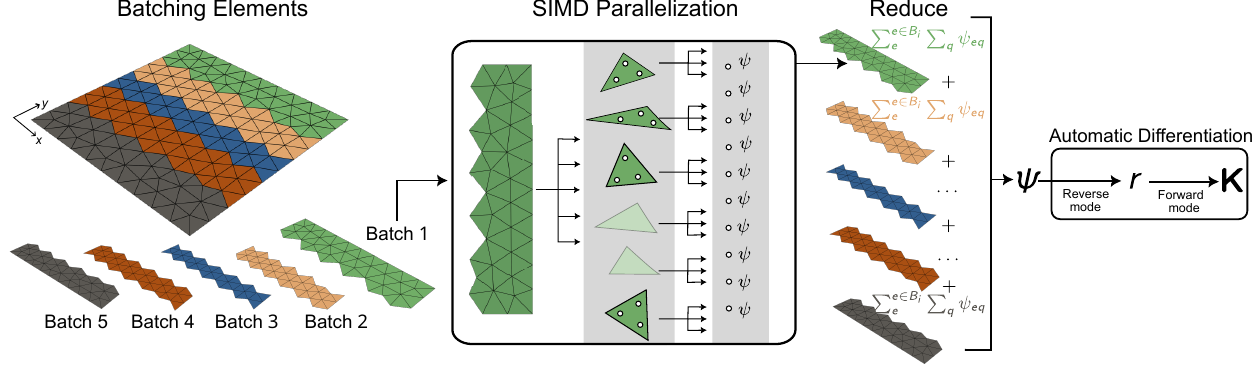}
	\caption{Schematic of the batched monolithic engine. To mitigate memory bottlenecks on GPUs, the global mesh is partitioned into contiguous element batches. Each batch is processed via SIMD-parallelized kernels that evaluate the local energy density $\psi$ simultaneously at each quadrature point from all elements in that batch and reduce them to a potential energy for that batch. Energy contribution from each batch is further reduced to get the global potential energy $\Psi_h$. The framework then recovers the global residual vector $\boldsymbol{r}$ (via reverse-mode AD) and the stiffness operator $\mathbf{K}$ (via forward-mode AD) directly from this scalar graph, bypassing the legacy element-level matrix assembly entirely.}
	\label{fig:batched-parallelization}
\end{figure}

In the proposed framework, we handle the transition from the mathematical formulation to hardware execution through batched computation. Instead of iterating over elements individually, we group elements into contiguous batches and process them collectively, see \Cref{fig:batched-parallelization}.  By unrolling the element loop into smaller groups, as shown in \Cref{alg:batched}, we ensure the peak memory footprint depends on the batch size rather than the total mesh size, preventing the materialization of the full computational graph for large problems. This execution strategy aligns naturally with the SIMD (Single Instruction, Multiple Data) model of modern processors and shifts the finite-element computation toward a compute-bound regime. As a result, the framework maintains high hardware occupancy through internal vectorization while exposing a high degree of parallelism. As shown in Sec.~\ref{sec:performance}, this formulation sustains constant throughput and linear $\mathcal{O}(N)$ scaling on GPUs, even for simulations with millions of degrees of freedom.

\begin{algorithm}[!t]
\caption{Batched Energy Computation}
\label{alg:batched}
\begin{algorithmic}[1]
\Require Global vector $\boldsymbol{u}_h\in\mathbb{R}^{N_{\mathrm{dofs}}}$, geometric dimension $d$, degrees of freedom per node $m$.
\Ensure Total potential $\Psi_h \in\mathbb{R}$.
\State Let $a \in \{1,\dots, n_{\mathrm{en}}\}$ denote the local node index per element.
\State Let $q \in \{1,\dots,N_q\}$ denote the local quadrature point index per element.
\State Reshape to nodal representation: $\boldsymbol{u}_h\rightarrow \boldsymbol{u}\in\mathbb{R}^{N_{\mathrm{nodes}}\times m}$.
\State \textbf{Partition} elements $\{1,\dots,N_{\mathrm{el}}\}$ into batches $\{B^k\}$ of size $N_B$.
\State $\Psi_h \gets 0$
\For{\textbf{each} batch $B^k$} \Comment{e.g., via \texttt{jax.lax.map}}
\State \textbf{Gather} geometry factors $\det(\mathbf{J}^k(\xi_q))$. 
  \Comment{Shape: [$N_B$, $N_q$]}
  \State \textbf{Gather} element nodal values $\boldsymbol{u}^k$. \Comment{Shape: [$N_B$, $n_{\mathrm{en}}$, $m$]}
  \State \textbf{Gather} element basis gradients $\nabla \mathcal{N}^k(\xi_q)$ \Comment{Shape: [$N_B$, $N_q$, $n_{\mathrm{en}}$, $d$]}
  
  \State \textbf{Compute} field gradients at quadrature points:
  \[
    \nabla \boldsymbol{u}^k(\xi_q)
    =
    \sum_{a=1}^{n_{\mathrm{en}}}
    \boldsymbol{u}_a^k \otimes \nabla \mathcal{N}_a^k(\xi_q).
  \]
  \Comment{Shape: [$N_B$, $N_q$, $m$, $d$]}
  \State \textbf{Evaluate} strain energy density $\psi^{k}(\xi_q)=\psi(\nabla \boldsymbol{u}^k(\xi_q))$.
  \Comment{Shape: [$N_B$, $N_q$]}
  \State \textbf{Integrate} element energies:
  \[
    \boldsymbol{\Psi}^{k}
    =
    \sum_{q=1}^{N_q}
    \psi^{k}(\xi_q)\, w_q\, \det(\mathbf{J}^k(\xi_q)).
  \]
  \Comment{Shape: [$N_B$]}
  \State \textbf{Accumulate} $\Psi_h \gets \Psi_h + \sum \boldsymbol{\Psi}^{k}$. \Comment{Scalar}
\EndFor
\end{algorithmic}
\end{algorithm}

We apply AD directly to the global scalar potential $\Psi_h(\boldsymbol{u}_h)$ to derive the residual vector and the tangent stiffness matrix. By differentiating the scalar energy functional, we ensure that all computational operators remain consistent with the fundamental variational formulation, see \cref{eq:fundamental-equation}. Operating at the level of the global energy allows us to express the full physics in a single computational graph, independent of any element-wise assembly logic. To make this approach computationally efficient, we select the AD mode based on the mapping dimensionality at each stage of the pipeline, thereby balancing computational cost and memory usage.

To compute the residual vector $\boldsymbol{r}$, we use reverse-mode AD. The global energy functional maps a high-dimensional input field $\boldsymbol{u}$ to a single scalar value, making reverse-mode AD the most efficient choice for this $N_{\text{dofs}} \rightarrow 1$ mapping~\cite{pundir_simplifying_2025}. Using this mode, we obtain the full residual vector at a computational cost that is only a small multiple of the cost of a single energy evaluation, independent of the total number of degrees of freedom.

To compute the tangent stiffness matrix $\mathbf{K}$, we need to evaluate second-order derivatives of the global energy functional. A direct application of AD at this level typically produces a dense Hessian, whose storage scales quadratically with the number of degrees of freedom. For large-scale problems, this dense representation quickly becomes infeasible due to high memory consumption. Consequently, a naive global differentiation of the energy functional does not provide a practical path towards scalable stiffness evaluation.

To address this challenge, we use two complementary strategies that operate directly on the global energy functional. When only the action of the stiffness matrix is required, for example, in iterative solvers or time integration schemes, we rely on matrix-free Jacobian–vector products, as described in \Cref{sec:jvp}. In contrast, when an explicit stiffness matrix is needed, such as for eigenvalue analyses or saddle-point problems, we recover a sparse representation of the tangent operator using graph-coloring-based differentiation, as discussed in \Cref{sec:coloring}. Together, these strategies allow the framework to remain scalable while supporting a broad range of problem classes.

\subsection{Matrix-Free Jacobian-Vector Products (JVP)}
\label{sec:jvp}
For matrix-free solvers, we compute the action of the tangent stiffness matrix expressed as Jacobian-Vector Products (JVP), without explicitly forming the matrix. We obtain JVPs by applying forward-mode AD to the residual function,
\begin{equation} \label{eq:single-jvp}
\delta \boldsymbol{r} = \underbrace{\frac{\partial \boldsymbol{r}}{\partial \boldsymbol{u}}}_{\mathbf{K}} \delta \boldsymbol{u}  \Rightarrow \delta \boldsymbol{r} = \text{JVP}(\delta \boldsymbol{u})~.
\end{equation}
This operation evaluates the directional derivative of the residual in the direction $\delta \boldsymbol{u}$ and naturally from the monolithic formulation, in which differentiation is applied at the level of the global functional. 

Matrix-free JVPs are well-suited for large-scale simulations in which iterative solvers or time integration schemes require repeated applications of the stiffness operator but do not require explicit access to the stiffness matrix. By avoiding the materialization of $\mathbf{K}$, this approach substantially reduces memory usage, and the cost of a single JVP evaluation remains comparable to that of computing the residual. Consequently, JVP evaluations scale linearly with the number of degrees of freedom.

\begin{figure}[!t]
	\centering
	\includegraphics[width=\linewidth]{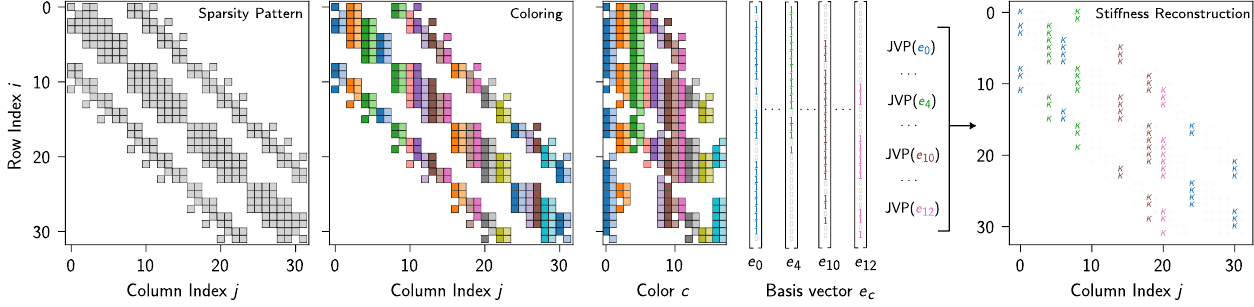}
	\caption{Schematic view of sparse differentiation via coloring algorithm. (from left-to-right) Foremost, we construct the sparsity pattern of a tangent stiffness matrix using mesh connectivity. The gray boxes represent the non-zero entries, and white represents zero entries. Next, we color the non-zero entries, \ie{} the gray boxes, using the distance-2 greedy algorithm. Non-zero entries in each color group are independent of each other and do not share any degrees of freedom. We group the entries color-wise into a compressed form with size $N_\mathrm{dofs} \times c$. Each column represents a basis vector (or composite perturbation vector) with 1 at non-zero entries and 0 in all other rows. Finally, applying a single Jacobian-vector product to each basis vector allows us to construct a compressed Jacobian matrix $\mathbf{J}_\text{comp}$ of size $N_\mathrm{dofs} \times c$ which is then used to reconstruct the sparse tangent stiffness matrix $\mathbf{K}$. }
	\label{fig:coloring-comparison}
\end{figure}

\subsection{Sparse differentiation using coloring} \label{sec:coloring}

For problems or solvers that rely on an explicit representation of the tangent stiffness matrix, we construct one using sparse differentiation~\cite{montoison_revisiting_2025, gebremedhin_what_2005} with \textit{matrix coloring}~\cite{curtis_estimation_1974}. In principle, this can be done by evaluating the JVPs of the residual $\boldsymbol{r}$ with respect to a set of basis vectors. Evaluating a JVP with a canonical basis vector $\boldsymbol{e}_i$ (1 at $i^{th}$ column, 0 elsewhere) yields the $i^\text{th}$ column of $\mathbf{K}$. Repeating this procedure for all $n$ degrees of freedom would recover the full matrix, but at a prohibitive cost of $n \times t_{\text{JVP}}$.

We reduce this cost by exploiting the sparsity structure of $\mathbf{K}$, as illustrated in \Cref{fig:coloring-comparison} and outlined in \Cref{alg:sparse-differentiation}. Using distance-based graph coloring of the finite-element connectivity, we partition the degrees of freedom into $c$ independent sets such that no two degrees of freedom within the same set contribute to the same row of the stiffness matrix. Each color, therefore, defines a composite perturbation vector, and a single JVP evaluated with this vector recovers multiple non-overlapping columns of $\mathbf{K}$ simultaneously.

The total cost of reconstructing the stiffness matrix is thus proportional to the number of colors, $c$, rather than the total number of degrees of freedom. For standard finite-element discretizations, $c$ is bounded by a constant that depends only on the element type and local connectivity and is independent of the global mesh size (see Appendix \ref{appendix:coloring} for more details). In practice, this results in a small number of JVP evaluations, even for meshes with millions of degrees of freedom. Because all degrees of freedom within a color are independent, each color pass can be evaluated in parallel without memory conflicts, allowing the compiler to exploit hardware parallelism efficiently. Furthermore, each pass to JVP can be evaluated in parallel because there is no memory conflict among the groups.

\begin{algorithm}
\caption{Sparse differentiation using coloring}
    \begin{algorithmic}[1]
        \Require Evaluation point $\boldsymbol{u}$, residual function $\boldsymbol{r}$, DOF color mapping $\boldsymbol{c}$, number of colors $n_{\text{colors}}$.
        \Ensure Compressed Jacobian $\mathbf{J}_{\text{comp}}$.
        \State \textbf{Initialize} color IDs: $\{0, 1, \dots, n_{\text{colors}}-1\}$.
        \For{\textbf{each} $c \in \text{color IDs}$} \Comment{Implemented via \texttt{jax.lax.scan}}
            \State \textbf{Generate} seed vector $\boldsymbol{e}$: $e_j = 1$ if $c_j = c$, else $0$.
            \State \textbf{Compute} JVP: $\boldsymbol{y} = \nabla \boldsymbol{r}(\boldsymbol{u}) \cdot \boldsymbol{e}$. \Comment{Directional derivative}
            \State \textbf{Store} $\boldsymbol{y}$ as column $c$ of $\mathbf{J}_{\text{comp}}$.
        \EndFor
        \State \Return $\mathbf{J}_{\text{comp}}$ with shape $[N_{\text{dofs}}, n_{\text{colors}}]$.
    \end{algorithmic}

\vspace{0.5em}

    \begin{algorithmic}[1]
        \Require Compressed Jacobian $\mathbf{J}_{\text{comp}}$, row pointers, column indices, color mapping $\boldsymbol{c}$.
        \Ensure Sparse stiffness matrix $\mathbf{K}$ in sparse format.
        \State \textbf{Expand} row pointers to generate row indices $i$ for every non-zero entry from sparsity pattern.
        \State \textbf{Retrieve} column indices $j$ from sparsity pattern.
        \State \textbf{Map} column indices to colors: $c = \boldsymbol{c}[j]$.
        \State \textbf{Extract} values: $K_{ij} = \mathbf{J}_{\text{comp}}[i, c]$. 
        \State \textbf{Construct} sparse matrix $\mathbf{K}$ from indices $(i, j)$ and values $K_{ij}$.
        \State \Return $\mathbf{K}$.
    \end{algorithmic}
    \label{alg:sparse-differentiation}
\end{algorithm}

The proposed framework, as described in this section, defines a monolithic, energy-centric approach in which physical models are expressed entirely in terms of a global scalar potential and differentiated systematically using automatic differentiation. By combining batched evaluation of the energy functional with global differentiation strategies, the framework avoids traditional element-wise assembly while remaining compatible with modern hardware architectures. The availability of both matrix-free Jacobian–vector products and sparse stiffness recovery provides flexibility in how second-order derivatives are accessed, allowing the same formulation to support a wide range of problem classes without modifying the underlying physics description.

\section{Versatility of the proposed framework}

The proposed framework is designed to isolate the stiffness construction from the specific physical model, allowing a wide range of problems to be expressed within a single computational architecture. This design choice enables the same formulation to be applied across diverse classes of problems, as demonstrated in the examples presented in \Cref{sec:examples}. Although the framework is formulated in terms of a potential energy functional, it is not restricted to conservative systems. As discussed in \Cref{sec:non-variational}, the architecture naturally extends to non-variational problems, such as fluid transport or non-associative plasticity. Furthermore, the same monolithic structure facilitates the integration of data-driven models without modifying the core formulation, as described in \Cref{sec:future_proof}.

\subsection{Generalization to non-variational systems}\label{sec:non-variational}

Our architecture handles the non-variational problems by extending the energy-centric approach to the principle of virtual work, as we show through an example in \cref{sec:ex:non-variational}. We introduce a scalar virtual work functional $\mathcal{W}(\boldsymbol{u},\boldsymbol{v})$, which represents the work done by internal and external forces against a test function $\boldsymbol{v}$: 
\begin{equation} 
\mathcal{W}(\boldsymbol{u}, \boldsymbol{v}) = \int_{\Omega} \boldsymbol{\sigma}(\boldsymbol{u}) : \nabla \boldsymbol{v}  d\Omega - \int_{\Omega} \boldsymbol{f} \cdot \boldsymbol{v}  \mathrm{d\Omega} 
\end{equation} 
The discretized form follows the same logic as the energy formulation by summing scalar contributions over quadrature points: 
\begin{equation} \mathcal{W}_h(\boldsymbol{u}_h, \boldsymbol{v}_h) = \sum_{e=1}^{N_\mathrm{el}} \sum_{q=1}^{N_\mathrm{q}}\!\left( \boldsymbol{\sigma}(\boldsymbol{u}^e_h (\xi_q))\right) : \nabla \boldsymbol{v}_h^{e}(\xi_q)\,  w_q\, \det\!\left(\mathbf{J}^{e}(\xi_q)\right) 
\end{equation} 
Our batched computation algorithm works exactly the same way for these problems. The framework simply aggregates the scalar product $\sigma:\nabla \boldsymbol{v}$ instead of the energy density $\Psi_h$. We then recover the global residual and tangent stiffness by differentiating this functional with respect to the test and trial functions: 
\begin{equation} \boldsymbol{r}(\boldsymbol{u}) := \nabla_{\boldsymbol{v}} \mathcal{W}(\boldsymbol{u}, \boldsymbol{v}) \Big|_{\boldsymbol{v}=\boldsymbol{0}}, \quad \mathbf{K}(\boldsymbol{u}) := \nabla_{\boldsymbol{u}} \boldsymbol{r}(\boldsymbol{u}) = \nabla^2_{\boldsymbol{u}\boldsymbol{v}} \mathcal{W}(\boldsymbol{u}, \boldsymbol{v}) \bigg|_{\boldsymbol{v}=\boldsymbol{0}}
\end{equation} 
This approach ensures that the monolithic architecture supports operator fusion and GPU acceleration for non-symmetric or path-dependent physics.

\subsection{Future-proofing with implicit AI coupling} \label{sec:future_proof}

The monolithic architecture also simplifies the integration of data-driven components. Traditional approaches~\cite{latyshev_expressing_2025, alheit_commet_2026} often treat machine learning models as external plugins, which requires batching local updates within a legacy finite-element loop. These methods rely on custom bridge code and manual derivations to communicate between the model and the solver. In contrast, our framework treats data-driven models as core physical entities, as demonstrated on two applications in \Cref{sec:integrate_data_driven}. Whether we integrate local Neural Constitutive Models (NCM)~\cite{thakolkaran_nn-euclid_2022, joshi_bayesian-euclid_2022,han_learning_2025} at the quadrature level, or domain-scale Neural Operator Element Methods (NOEM)~\cite{wang_accelerating_2025}, the automatic differentiation treats the neural network as an additive term in the global scalar sum. This implicit AI coupling provides a unified structure that remains scalable for high-performance computing while staying flexible enough for the next generation of AI-augmented mechanics.

\section{\tatva: A monolithic differentiable library}\label{sec:tatva}

While the proposed energy-centric framework is library-agnostic, we have implemented these concepts into an open-source Python library called \tatva~\cite{tatva_library_2026}. Built on the JAX ecosystem~\cite{bradbury_jax_2018, kidger_equinox_2021, gerdes_mathisgerdesjax-autovmap_2025}, \tatva provides a high-level API for variational mechanics while maintaining the performance of a low-level GPU-accelerated framework. The primary goal is to allow researchers to define complex physical potentials with minimal code. The library centers on the \texttt{Operator} class, which manages the relationship between mesh topology, element formulation, and the global state.

We summarize the key functionalities below:

\begin{itemize}
\item \textbf{Geometric Context:} The \texttt{Operator(mesh, element, batchsize)} object initializes the element connectivity, shape functions, and quadrature rules. It sets the batch size for SIMD operations such as integration and gradient calculations.

\item \textbf{SIMD Batched Operations:} The \texttt{Operator} object provides high-level functions such as \texttt{op.grad} and \texttt{op.eval}, which use \texttt{jax.lax.map} to batch computations across the mesh. \texttt{op.grad} evaluates nodal gradients at all quadrature points, and \texttt{op.eval} interpolates nodal values. We use \texttt{jax.lax.map} instead of \texttt{jax.vmap} to preserve memory scalability: \texttt{vmap} materializes the full differentiable graph across all elements, quickly exhausting memory for large problems. In contrast, \texttt{jax.lax.map} iterates over element batches, bounding peak memory while still allowing the XLA compiler~\cite{noauthor_openxlaxla_nodate} to fuse quadrature-level operations into optimized kernels.

\item \textbf{Variational Integration:} The \texttt{op.integrate} function performs the final reduction. It takes the local energy density contributions or any scalar value at all quadrature points and applies the appropriate Jacobian and quadrature weights to produce the global scalar potential $\Psi_h$. 
\end{itemize}

The design of the library is a direct realization of the ``gather-evaluate-reduce'' pipeline described in \Cref{alg:batched}. A typical implementation of the total strain energy translates the pseudo-algorithm into executable code line-by-line:

\begin{lstlisting}
op = Operator(mesh, element_type) # define the operator for a mesh and a element type

def total_strain_energy(u_flat):
    u = u_flat.reshape(-1, n_dofs)                # reshaping
    u_grad = op.grad(u)                           # Compute Gradients in batches
    energy_density = strain_energy_density(u_grad)  # Quadrature-level Evaluation
    return op.integrate(energy_density)            # Reduction via Batched Integration

u_flat = ... # 1D array of length equal to total DoFs
total_energy = total_strain_energy(u_flat)
\end{lstlisting}

In this energy function, we pass degrees of freedom as a flattened array. This choice ensures that when the library automatically differentiates the function, it generates a residual vector and a stiffness matrix with the correct dimensions ($N_\text{dofs}$  and $N_\text{dofs} \times N_\text{dofs}$) respectively.

Because \tatva uses \texttt{jax.numpy} for its data structures, handling and slicing data is as simple as using standard NumPy. This compatibility allows researchers to express complex mathematical expressions, conditional logic, or loops easily to compute energy density from quadrature-level values. This structure makes it possible to represent sophisticated physics with minimal coding effort while keeping the implementation closely aligned with the original mathematical theory.

Furthermore, the design ensures that every derivative computed by the library, whether it is the residual $\boldsymbol{r}$ or the stiffness matrix $\mathbf{K}$, is derived from the same monolithic computational graph. This consistency allows \tatva to switch between colored sparse assembly and the matrix-free JVP approach without any modifications to the user's physics code. From the energy function defined above, we compute the residual function using \texttt{jax.jacrev} directly:

\begin{lstlisting}
compute_residual = jax.jacrev(total_energy)
residual = compute_residual(u_flat)
\end{lstlisting}

One can then utilize this residual function to compute the action of the stiffness matrix (see Eq. \ref{eq:single-jvp}) via \texttt{jax.jvp}:

\begin{lstlisting}
tangent = jax.jvp(compute_residual, (u_flat), (du_flat))[1]
\end{lstlisting}

To construct the sparse stiffness matrix, \tatva provides functionalities to generate sparsity patterns and colors using a greedy distance-2 algorithm. These colors are then used to compute the sparse matrix as described in \Cref{alg:sparse-differentiation}:

\begin{lstlisting}
from tatva import sparse # importing sparse differentiation module

# setup sparsity and coloring
sparsity_pattern = sparse.create_sparsity_pattern(mesh, n_dofs_per_node)
colors = sparse.distance2_colors(sparsity_pattern) # vector with non-zero entries colored

# materialize the sparse stiffness matrix
K_sparse = sparse.jacfwd(compute_residual, sparsity_pattern, colors) # Algorithm 2 
\end{lstlisting}

This architecture provides a truly universal infrastructure for differentiable mechanics, enabling the user to pivot between robustness and scalability at the solver level without ever revisiting the physics implementation.

\section{Performance analysis} \label{sec:performance}

In this section, we analyze the computational efficiency and scalability of the proposed framework to construct tangent operators on modern hardware. We benchmark the framework against traditional assembly-based methods to highlight the advantages of our graph coloring strategy and matrix-free JVPs. Our analysis uses two setups: a 2D elasticity problem with linear triangular elements and a 3D elasticity problem with tetrahedral elements. By varying the mesh resolution, we assess scaling behavior across a wide range of degrees of freedom (DoFs). We measure the execution time and computational throughput, measured in million DoFs per second (million DoFs/s), for two operations: matrix-free JVP evaluation and sparse assembly via graph coloring. All benchmarks were performed on an NVIDIA A100 GPU.

\begin{figure}[t]
	\centering
	\includegraphics[width=\linewidth]{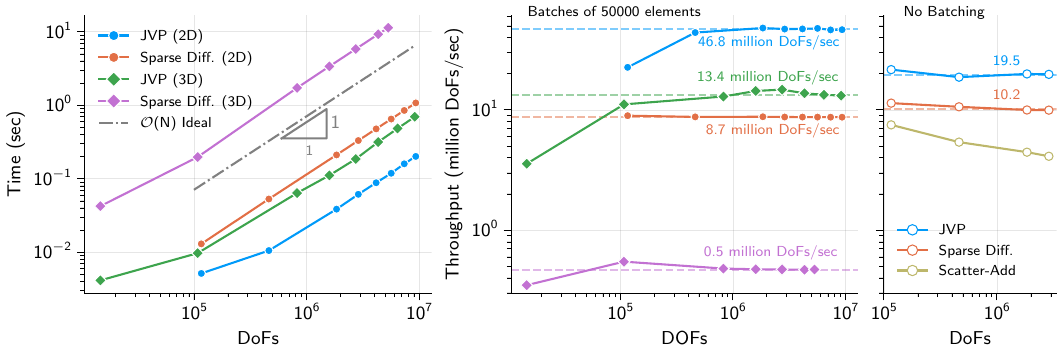}
	\caption{Computational performance and scaling of the energy-centric framework. (left) The \textit{log-log} shows execution time for a computation with batch size of 50,000 elements. The unit slope for all methods confirms strict $\mathcal{O}(N)$ scaling. (middle) Assembly throughput (Million DoFs/sec) for the same computations, demonstrating optimal GPU utilization across three orders of magnitude. (right) Comparison to a scatter-add assembly approach similar to~\cite{xue_jax-fem_2023, ichihara_naruki-ichiharafeax_2026} for a 2D computation without batching for fair comparison. The throughput for scatter-add assembly is significantly lower, no longer constant, and constantly drops with increasing DoFs, clearly a consequence of $\mathcal{O}(N^2)$ scaling. }
	\label{fig:performance}
\end{figure}

\subsection{Scaling and throughput analysis}

We evaluate the scalability and efficiency of the proposed framework by analyzing both execution time and computational throughput with increasing DoFs up to 10 million DoFs. For both matrix-free JVP evaluation and sparse differentiation via graph coloring, the results exhibit linear $\mathcal{O}(N)$ scaling over more than three orders of magnitude in problem size (see \Cref{fig:performance}-left), confirming that the proposed architecture avoids algorithmic overheads that grow with system size.

Expressed in terms of throughput (see \Cref{fig:performance}-center), the results demonstrate sustained GPU efficiency across all mesh resolutions. The matrix-free formulation achieves peak throughputs of approximately 47 million DoFs/s in 2D and 13 million DoFs/s in 3D, with throughput remaining essentially constant as the problem size increases. Sparse differentiation via graph coloring exhibits lower throughput due to the multiple JVP evaluations required per color, but remains stable at approximately 8.7 million DoFs/s in 2D and 0.5 million DoFs/s in 3D.

For the 3D sparse differentiation case, the available 40~GB of GPU memory limited the maximum problem size to approximately 5.1 million DoFs. This limitation highlights the importance of memory-efficient differentiation strategies. More advanced coloring techniques~\cite{montoison_revisiting_2025} could further reduce the number of required colors and extend this limit.

\subsection{Comparison with traditional assembly}

We compare the throughput of our approach with the traditional scatter-add assembly process (see \Cref{fig:performance}-right). For a fair comparison, we do not use any batching in either the reference approach or our proposed framework. The results show two significant advantages of the proposed energy-centric framework. First, both our matrix-free JVP and sparse differentiation methods achieve higher throughput than the traditional scatter-add approach. Second, unlike the constant throughput of our proposed approach, the throughput of scatter-add assembly constantly decreases with increasing number of DoFs. This performance drop is a clear indicator of the memory access bottlenecks and atomic contention inherent in element-centric assembly. In contrast, our proposed approach of evaluating the tangent operators maintains stable throughput, demonstrating that the monolithic design avoids these traditional limitations and scales effectively for large-scale engineering problems.

\section{Examples} \label{sec:examples}

To illustrate the versatility of the proposed framework, we apply the energy-centric formulation to a diverse set of complex physical problems. These examples demonstrate that a wide range of governing equations, constraints, and couplings can be expressed within a single unified implementation without intrusive modifications to the software stack. In all cases, the formulations rely only on three basic operations, \texttt{op.grad}, \texttt{op.eval}, and \texttt{op.integrate}, introduced in \Cref{sec:tatva}, highlighting the simplicity and generality of the approach.

We begin in \cref{sec:ex1:bcs} with (linear) elasticity and essential boundary conditions to validate our framework and establish a baseline for subsequent, more complex examples. In \cref{sec:mpc}, we enforce multi-point constraints using techniques such as condensation, Lagrange multipliers, and the penalty method. In \cref{sec:mixed}, we demonstrate mixed-dimensional coupling, including embedded fibers and cohesive interfaces. In \cref{sec:ex:non-variational}, we treat a non-variational problem formulated through the principle of virtual work. Finally, in \cref{sec:integrate_data_driven}, we integrate data-driven models using neural networks in the same energy-centric structure.

For readability, the code listings focus on the essential components of each formulation rather than complete scripts. Fully executable versions of all examples are available in the \tatva documentation~\cite{tatva_docs_2026}.

\subsection{Handling elasticity and boundary conditions}\label{sec:ex1:bcs}

\begin{figure}[!t]
    \centering
    \includegraphics[width=\linewidth]{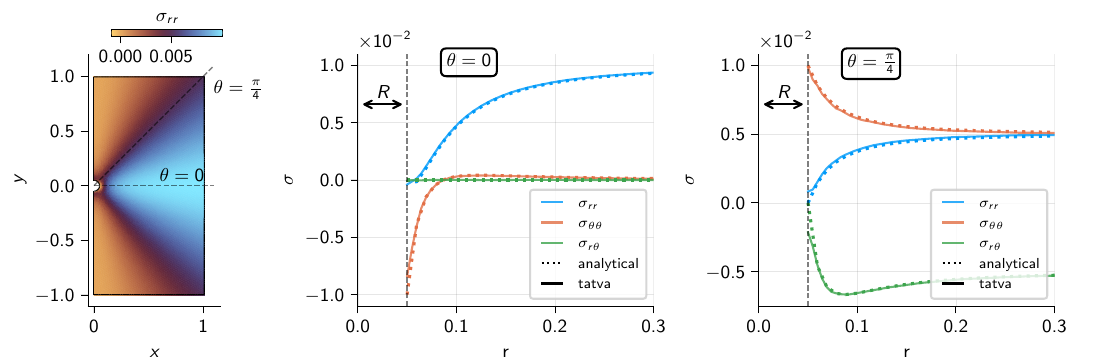}
    \caption{
    Plate with a circular hole of radius $R=0.05$ subjected to horizontal uniaxial tension $t_x=0.01$. 
    {(left)} Distribution of the radial stress component $\sigma_{rr}$ in the half-domain model with
    symmetry boundary conditions imposed on the left edge. Dashed lines indicate radial paths corresponding to $\theta=0^\circ$ (loading direction) and $\theta=45^\circ$. 
    {(middle/right)} Radial stress components $\sigma_{rr}$, $\sigma_{\theta\theta}$, and $\sigma_{r\theta}$ along $\theta=0^\circ$, and $\theta=45^\circ$, respectively, plotted as functions of the radial coordinate $r$. The solid lines represent the results obtained using the framework presented here; these results are compared with Kirsch’s analytical solution for an infinite plate with a circular hole (dashed line). 
    }
    \label{fig:ex_linear_elastic}
\end{figure}

We first consider a linear elasticity problem to illustrate how standard boundary-value problems are expressed within the proposed framework. The example consists of a plate with a circular hole at its center being subjected to a uniaxial tension, as shown in \cref{fig:ex_linear_elastic}-left. Owing to symmetry, we model only one half of the domain, with symmetry boundary conditions imposed along the left
edge. A prescribed traction is applied on the right boundary of the plate.

The problem can be formulated as the minimization of the total potential energy subject to kinematic constraints,
\begin{gather}
    \Psi(\boldsymbol{u}) 
    = \int_{\Omega} \psi_\varepsilon\!\left(\boldsymbol{\varepsilon}(\boldsymbol{u})\right)\,\mathrm{d}\Omega
    - \int_{\mathrm{S}_t} \boldsymbol{t} \cdot \boldsymbol{u}\,\mathrm{d}\Gamma ~, \\
    \boldsymbol{u}_x = \mathbf{0} \quad \text{on } \mathrm{S}_{D} ~,  \notag
\end{gather}
where $\psi_\varepsilon$ denotes the strain energy density and $\boldsymbol{t}$ the applied traction. 

In the proposed framework, this formulation translates directly into code by defining the strain energy density and assembling the total potential energy through the core operators as follows:
\begin{lstlisting}
def strain_energy_density(u_grad: Array) -> Array:    # linear elastic material
     eps = 0.5 * (u_grad + u_grad.T)
     sig = 2 * mu * eps + lmbda * jnp.trace(eps) * jnp.eye(2)
     return 0.5 * jnp.einsum("ij,ij->", sig, eps)

def total_energy(u_flat: Array, traction: float) -> float:
    u = reshape(u_flat)
    u_grad = op.grad(u)
    psi = strain_energy_density(u_grad)
    # op_right_edge is an Operator with the corresponding line elements
    line_integrand = op_right_edge.eval(traction * u)  
    return op.integrate(psi) - op_right_edge.integrate(line_integrand)
\end{lstlisting}

Dirichlet boundary conditions are enforced by eliminating the prescribed degrees of freedom at the level of the energy functional. Instead of modifying an assembled linear system, \eg, via row/column operations or direct elimination of constrained DoFs~\cite{belytschko_finite_2008}, we define a reduced functional that depends only on the free degrees of freedom. The reduced displacement vector is lifted to the full field by inserting the prescribed values at constrained locations before evaluating the total energy. 
AD of this wrapped functional yields the consistent residual and tangent, such that the constraint handling can be carried out implicitly at the variational level:

\begin{lstlisting}
free_dofs_indices = jnp.where(mesh.coords[:, 0] != 0.0)[0]

def total_energy_reduced(u_free: Array, traction: float) -> float:
    u_full = jnp.zeros(n_dofs).at[free_dofs_indices].set(u_free)
    return total_energy(u_full, traction)

residual_condensed = jax.jacrev(total_energy_reduced)
K_condensed = jax.jacfwd(residual_condensed)
\end{lstlisting}

To confirm the accuracy of our model, the obtained numerical results are compared with Kirsch's analytical solution for an infinite plate with a circular hole under remote tension~\cite{kirsch1898theorie}. Excellent agreement is observed for all stress components, as shown in~\cref{fig:ex_linear_elastic}. 

To highlight the modularity of the framework, we note that the linear elastic model could be replaced by a finite-strain hyperelastic formulation by modifying only the strain energy density $\psi$. All other components, including discretization, quadrature, boundary-condition handling, and solution procedure, would remain unchanged. In practice, this requires defining a new mapping from the displacement gradient to a scalar energy density,
\begin{equation*}
\psi_\varepsilon : \nabla \boldsymbol{u} \mapsto \mathbb{R},
\end{equation*}
which is evaluated at the quadrature points and assembled into the global energy:

\begin{lstlisting}
def strain_energy_density(u_grad: Array) -> Array:   # hyperelastic material
    F = jnp.eye(2) + grad_u 
    C = F.T @ F   # computing Cauchy strain
    ... # do stuff
    return energy_density
\end{lstlisting}

The residual and consistent tangent operators are again obtained by AD of the assembled energy functional, without requiring any further changes to the implementation.

\begin{figure}
    \centering
    \includegraphics[width=\linewidth]{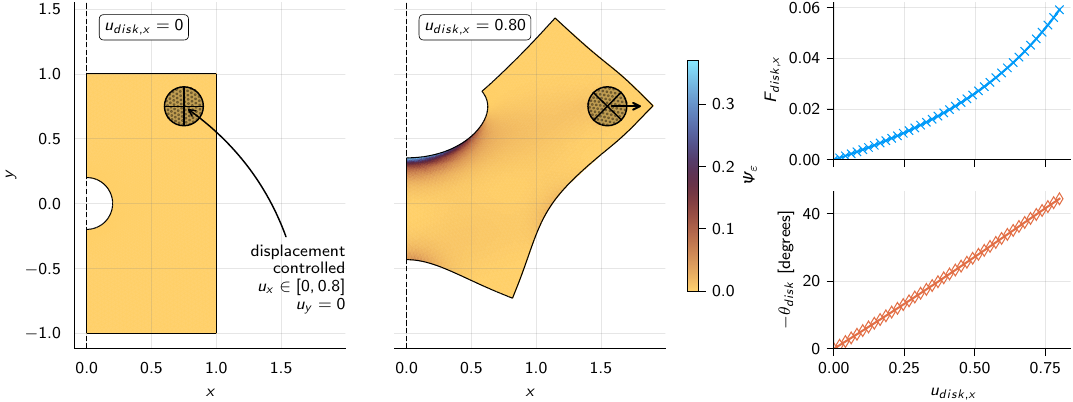}
    \caption{
    Rigid-body multi-point constraint (MPC) applied to a hyperelastic plate. 
    (left) Reference configuration with a rigid disk attached to the top-right corner. The disk motion is described by the generalized coordinates $(u_{\text{disk},x},u_{\text{disk},y},\theta_\text{disk})$, with the horizontal displacement $u_{\mathrm{disk},x}$ prescribed and $u_{\mathrm{disk},y}=0$. 
    Nodes within a radius $R_\mathrm{disk}=0.15$ of the disk center are constrained to follow rigid-body kinematics. 
    (middle) Deformed configuration at $u_{\text{disk},x}=0.8$, colored by the strain energy density $\psi_\varepsilon$. A symmetry boundary condition is imposed along the left edge (dashed line). 
    (right) Reaction force $F_{\text{disk},x}$ and rotation $\theta_{\text{disk}}$ of the rigid disk as functions of the imposed displacement $u_{\text{disk},x}$.
    }
    \label{fig:ex_nonlinear_mpc}
\end{figure}

\subsection{Multi-point constraints}\label{sec:mpc}
Multi-point constraints (MPCs) impose kinematic relations between the degrees of freedom associated with different nodes. Such constraints arise in a wide range of applications, including rigid-body motion, periodic boundary conditions for homogenization, and contact between deformable bodies. The proposed framework supports multiple strategies for enforcing MPCs, including static condensation, Lagrangian multiplier formulations, and penalty-based approaches, without requiring changes to the underlying differentiation pipeline.

\subsubsection{Rigid body constraints using condensation}\label{sec:ex2:rigid}
As a first illustration, we consider rigid-body motion enforced through static condensation, following the same principle used to impose essential boundary conditions in \cref{sec:ex1:bcs}. In this setting, we impose linear kinematic relations on a subset of nodes so that they follow a prescribed rigid-body displacement.

Here, we again consider a plate with a hole in its center; in this case, however, there is a discontinuity going from the center of the hole upwards to the top boundary (see \cref{fig:ex_nonlinear_mpc}) and the material is assumed to be hyperelastic. Owing to symmetry, we only model the right half of the domain. We attach a rigid disk of radius (centered at $x=0.75, y=0.75$) to the top-right part of the plate and displace it to induce deformation in the surrounding material. We describe the motion of the disk using the generalized degrees of freedom $(u_{\mathrm{disk},x}, u_{\mathrm{disk},y},\theta_{\mathrm{disk}})$, where $u_{\mathrm{disk},x}$ and $u_{\mathrm{disk},y}$ are prescribed, and $\theta_{\mathrm{disk}}$ is free to adapt such that the disk does not experience any moment, and hence $\theta$ is an unknown. 
All nodes located within a radius $R_\mathrm{disk}=0.15$ of the disk center follow the rigid-body kinematics of the disk through MPC relations. 

We construct the vector of unknowns by combining the nodal displacements and the disk degrees of freedom according to,
$$
\boldsymbol{z} = \begin{bmatrix}
    \boldsymbol{u} \\ \boldsymbol{u}_\text{disk}
\end{bmatrix} ~ .
$$
As in the case of essential boundary conditions, we enforce the constraints by defining a reduced energy functional that depends only on the free degrees of freedom. We reconstruct the full displacement field by lifting the reduced variables and applying the rigid-body kinematic relations before evaluating the energy:

\begin{lstlisting}
def apply_rigid_body_bc(u_full: Array, u_disk: Array) -> Array:
    ux_c, uy_c, theta_c = u_disk # get the 3 rigid components
    cos, sin = jnp.cos(theta_c), jnp.sin(theta_c)  # construct rotation matrix based on theta
    rotation_matrix = jnp.array([[cos, -sin], [sin, cos]])
    vecs = mesh.coords[disk_nodes] - disk_center 
    rotated_vecs = vecs @ rotation_matrix.T # rotate vectors to the disk center
    u_rotation = rotated_vecs - vecs        # get the displacement vector
    u_rigid = jnp.stack([ux_c, uy_c]) + u_rotation  # add to the rigid body translation
    u_full = u_full.at[disk_nodes_x].set(u_rigid[:, 0]) # set u in x direction
    u_full = u_full.at[disk_nodes_y].set(u_rigid[:, 1])  # set u in y direction
    return u_full

def total_energy_reduced(u_free: Array, u_disk_x: float) -> float:
    u_full = jnp.zeros((n_dofs,)).at[free_dofs].set(u_free)
    u_disk = u_full[-3:] # we extract the 3 rigid body DoFs
    u_disk = u_disk.at[0:2].set([u_disk_x, 0]) # set the constraints
    u_full = apply_rigid_body_bc(u_full, u_disk)
    return total_energy(u_full)
\end{lstlisting}

\Cref{fig:ex_nonlinear_mpc} shows the resulting deformed configuration, together with the force–displacement and rotation–displacement responses of the rigid disk.

\subsubsection{Periodic constraints for homogenization using Lagrange multipliers}\label{sec:lagrange}

This example demonstrates the implementation of periodic boundary conditions for the homogenization of a linear elastic composite using Lagrange multipliers. We consider a two-dimensional plane strain problem comprised of a skewed unit cell containing stiff circular inclusions of radius $R$, with elastic properties $(E_r,\nu_r)$, embedded in a compliant linear elastic matrix with properties $(E_m,\nu_m)$. The inclusions are arranged in a hexagonal pattern, for which homogenization theory predicts an isotropic effective macroscopic response. The example is inspired by~\cite{bleyer2018numericaltours}.

We impose periodic boundary conditions on the displacement fluctuation field using Lagrange multipliers. The multipliers are included in the global vector of unknowns $\boldsymbol{z} = (\boldsymbol{u}, \boldsymbol{\lambda})$, and the Lagrangian
functional reads
\begin{equation}
\mathcal{L}(\boldsymbol{u}, \boldsymbol{\lambda})
= \Psi(\boldsymbol{u}) + \boldsymbol{\lambda} \cdot \boldsymbol{g}(\boldsymbol{u})~,
\end{equation}
where $\Psi(\boldsymbol{u})$ is the total energy and $\boldsymbol{g}(\boldsymbol{u})$ collects the periodicity constraints enforcing equality of displacement fluctuations on opposing boundaries of the unit cell.

In the proposed framework, we express these constraints directly at the level of the global functional:
\begin{lstlisting}
def g(u: Array) -> Array:
    """Returns the constraint violations."""
    g_left_right = u[right] - u[left]  # right/left are paired boundary nodes
    g_bot_top = u[top] - u[bot]        # bot/top are paired boundary nodes
    return jnp.concatenate([g_left_right, g_bot_top])
    
def total_energy(z_flat: Array) -> float:
    u, lmbda = slice_and_reshape_z(z_flat)
    ...  # computing the strain energy density
    energy_density = strain_energy_density(u_grad) # from Sec 6.1
    return op.integrate(energy_density) + jnp.dot(lmbda, g(u))
\end{lstlisting}

\begin{figure}[!t]
    \centering
    \includegraphics[width=\linewidth]{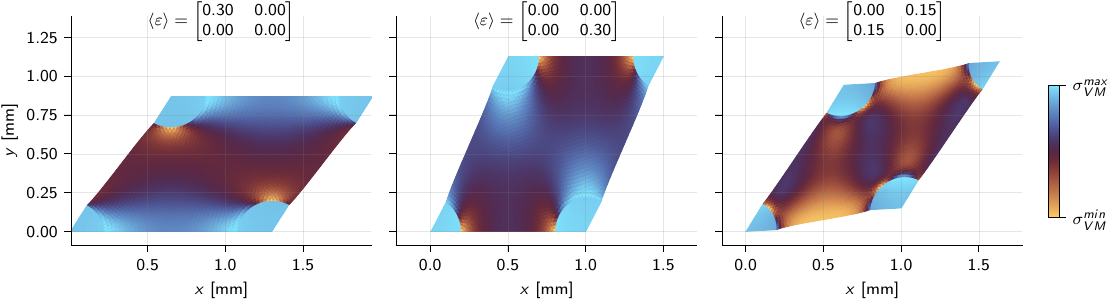}
    \caption{
        Periodic homogenization of linear elasticity with a skewed unit cell. The RVE consists of two phases; a softer matrix and a stiffer inclusion. Evaluations of the microscopic response for different applied average strain tensors $\langle\varepsilon\rangle$, corresponding to a pure stretch in $x$, $y$, and pure shear. This example is adopted from~\cite{bleyer2018numericaltours} with the same geometric and material parameters.
    }
    \label{fig:ex_homogenization}
\end{figure}

Stationarity of $\mathcal{L}$ with respect to both $\boldsymbol{u}$ and $\boldsymbol{\lambda}$ yields a saddle-point system, which we solve monolithically using a direct linear solver. We construct the corresponding sparse tangent operator using the graph-coloring-based differentiation approach described in \cref{sec:coloring}. Appendix~\ref{appendix:sparsity} details the procedure used to extract the sparsity pattern of the Lagrangian functional required for sparse differentiation. \Cref{fig:ex_homogenization} shows the deformed unit cell for three principal macroscopic strain states.

Beyond solving the microscopic equilibrium problem, the fully differentiable formulation enables the direct computation of the exact consistent homogenized constitutive tensor $\mathbb{C}_{\text{hom}}$. We define the macroscopic stress as the volume average of the microscopic stress obtained from the equilibrium solution corresponding to a prescribed macroscopic strain. Since the entire solve procedure is differentiable, we obtain the homogenized tangent operator by differentiating the macroscopic stress mapping with respect to the applied strain~\cite{pundir_simplifying_2025},
$$
\boldsymbol{\sigma}_{\text{hom}}(\boldsymbol{\varepsilon}_{\text{macro}})~.
$$
This approach yields the exact algorithmic homogenized tangent, including all microstructural equilibrium effects, without deriving analytical sensitivity expressions. In practice, we implement this procedure by defining a differentiable mapping from the applied macroscopic strain to the homogenized stress response:
\begin{lstlisting}
def solve(eps_macro_voigt: Array) -> Array:
    eps_macro = from_voigt(eps_macro_voigt) # reshape voigt notation into 2x2 tensor
    u_tilde_flat = jnp.linalg.solve(jac(v0, eps_macro), -residual(v0, eps_macro))
    u_tilde = u_tilde_flat.reshape(-1, n_nodes)
    u = eps_hat @ x + u_tilde
    sig = jnp.mean(compute_stress(u), axis=0)  # average stress response
    return to_voigt(sig)  # return average stress in voigt notation

C_homogenized = jax.jacfwd(solve)(jnp.ones(3))
\end{lstlisting}

\subsubsection{Contact constraints using penalty approach}
\label{sec:ex_contact}

\begin{figure}[!t]
    \centering
    \includegraphics[width=\linewidth]{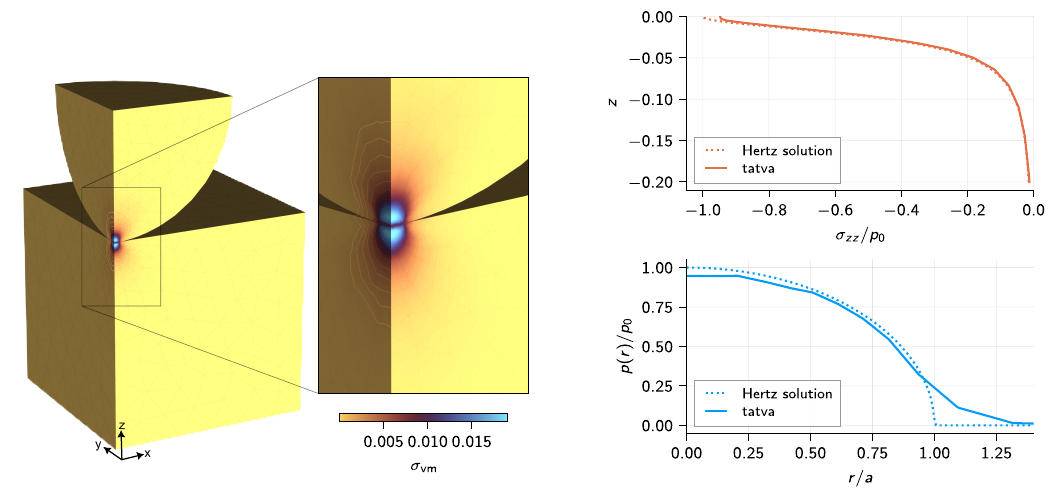}
        \caption{
        Three-dimensional Hertzian contact benchmark. The sphere radius is $R_\mathrm{Hertz}=0.6$, and both bodies are assumed to behave linear elastically with identical material parameters $E=1.0$ and $\nu=0.3$ (left) Quarter-symmetry model of a spherical indenter in frictionless contact with an elastic half-space; the inset shows a magnified view of the contact region. The deformed configuration is colored by the von Mises stress, $\sigma_{\mathrm{vm}}$. (right-top) Normalized subsurface stress component $\sigma_{zz}/p_0$ along the symmetry axis beneath the contact center, compared with the analytical Hertz solution. (right-bottom) Normalized contact pressure distribution $p(r)/p_0$ as a function of the radial coordinate $r/a$, where $a$ is the contact radius and $p_0$ the maximum contact pressure.
        }
    \label{fig:ex_contact}
\end{figure}

We next consider a three-dimensional contact problem to illustrate how the proposed framework handles geometric constraints that introduce non-local couplings between independently discretized bodies. We model a linear elastic sphere indenting a linear elastic cube (see ~\cref{fig:ex_contact}-left) and compare the numerical results with the analytical Hertz solution for frictionless contact~\cite{johnson1987contact}.

We discretize the two bodies independently and apply essential boundary conditions using condensation, as described in \cref{sec:ex1:bcs}. We fix the bottom surface of the cube in the $z$-direction and prescribe a vertical displacement on the top surface of the sphere octant. The interaction between the two bodies is introduced solely through the contact constraint.

Contact introduces geometric couplings that are not aligned with the mesh topology, as nodes on one surface interact with geometrically nearby but topologically unrelated regions of the opposing surface. In sparsity-based formulations, such long-range interactions must be explicitly reflected in the construction of the sparsity pattern. In contrast, the matrix-free formulation introduced here incorporates contact directly through an additional contribution to the global energy functional, without requiring any modification of the underlying algebraic structure.

We enforce the contact constraint using a penalty formulation by augmenting the total energy functional with a contact term $\Psi_c$ defined through a penalty density
$$
\psi_c =\frac{1}{2} \kappa~ g(\boldsymbol{u}_1, \boldsymbol{u}_2)^2~,
$$
where $g$ denotes the normal gap function and $\kappa$ is the penalty parameter. The total energy functional then becomes
\begin{equation}
\Psi(\boldsymbol{u}_1, \boldsymbol{u}_2) = \int_{\Omega_{\mathrm{sphere}}}\psi_\varepsilon(\boldsymbol{u}_1)~\mathrm{d\Omega} + \int_{\Omega_{\mathrm{cube}}}\psi_\varepsilon(\boldsymbol{u}_2)~\mathrm{d\Omega} +\int_{\mathrm{S}_c}\psi_c(\boldsymbol{u}_1, \boldsymbol{u}_2) ~\mathrm{dS} ~ .
\end{equation}

To evaluate the contact energy, we use a two-pass node-to-surface contact detection algorithm to determine the gap $g(\boldsymbol{u}_1, \boldsymbol{u}_2)$. In the first pass, we orthogonally project every surface node of the sphere onto the cube surface to compute the gap and then integrate the penalty density over the surface of the sphere. In the second pass, we compute the gap from every surface node of the cube to the sphere surface and then integrate the corresponding penalty energy over the surface of the cube. We define the total contact energy as the average of these two contributions:
\begin{lstlisting}
def compute_contact_energy(x1, x2):
    g1 = ... # gap for sphere surface nodes using node-to-surface contact algorithm
    g2 = ... # gap for cube surface nodes using node-to-surface contact algorithm
    contact_energy_1 = ... # integrate
    contact_energy_2 = ... # integrate
    return 0.5*(contact_energy_1 + contact_energy_2)

def total_energy(z_flat: Array) -> float:
    u1, u2 = slice_and_reshape_z(z_flat)
    # [...] Compute elastic energies for both bodies -> elastic_energy_1, elastic_energy_2
    # Compute contact energy based on two bodies' current configurations
    contact_energy = compute_contact_energy(x1 + u1, x2 + u2)
    return (op_1.integrate(elastic_energy_1) + op_2.integrate(elastic_energy_2)
        + contact_energy)
\end{lstlisting}

\Cref{fig:ex_contact} compares the numerical results with the analytical Hertz solution. Overall, we observe good agreement in both the subsurface stress distribution and the contact pressure profile. The deviations in the contact pressure outside the contact area are attributed to the relatively coarse discretization of the surface mesh.

We note that in this example, we deliberately employ a simple contact treatment to emphasize the generality of the energy-centric approach. More advanced contact methods, such as surface-to-surface or mortar methods~\cite{Popp2012}, can be incorporated in the same manner by modifying the contact energy term. Because AD automatically handles the resulting cross-coupling terms automatically, the framework substantially simplifies the implementation of contact problems involving complex geometric interactions.

\subsection{Mixed domain or mixed dimensional problems}\label{sec:mixed}

Mixed-domain or mixed-dimensional problems involve interactions between physical processes defined on regions of different topological dimensions, such as cohesive fracture along a 2D interface embedded in a 3D domain. Simulating such systems in standard finite-element frameworks can be challenging~\cite{daversin-catty_abstractions_2021}, as it requires managing multiple incompatible meshes, interpolating solution fields across non-conforming dimensions (\eg, from 3D bulk elements to 2D surface elements), and constructing complex assembly mappings to couple the different domains. 

In the proposed energy-centric framework, we address these challenges by treating the physics associated with each domain as an additive contribution to the global energy functional. Lower-dimensional contributions enter the formulation as additional energy terms defined on their respective geometries, while AD handles the resulting couplings consistently at the global level. This approach avoids the need for specialized assembly logic or explicit domain-to-domain mapping procedures. In the following examples, we illustrate how this formulation simplifies the implementation of mixed-dimensional topologies within a unified computational framework.

\begin{figure}[!t]
    \centering
    \includegraphics[width=\linewidth]{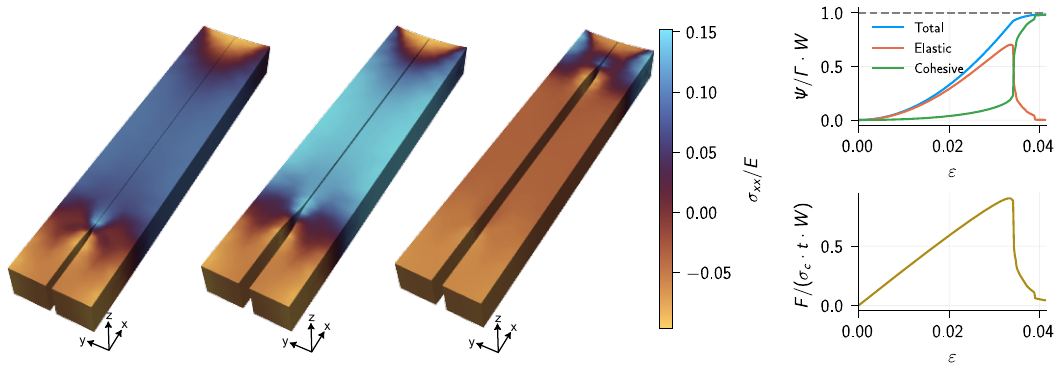}
    \caption{Cohesive fracture along an interface under Mode-I loading in the $y-$direction. A linear elastic body of dimensions $10L \times 2L \times L$ with $L=0.005$ and with an initial interfacial crack of length $L$ is loaded under displacement control. (left) Deformed configuration at three loading steps, illustrating crack propagation along the interface. Colors indicate the normal stress component $\sigma_{xx}$, normalized by the Young's modulus. (right-top) Evolution of the total strain energy, total fracture energy, and the total energy during loading. The total dissipated energy matches the theoretical value $\Gamma \times A_\text{crack}$ required to fully separate the interface. (right-bottom) Reaction force at the top edge as a function of the applied strain, with the sudden force drop indicating brittle fracture behavior.}
    \label{fig:cohesive-fracture}
\end{figure}

\subsubsection{Fracture using cohesive elements}

Cohesive fracture modeling typically relies on specialized interface elements~\cite{Ortiz1999, Vocialta2017a}, introducing a mixed-dimensional setting in which the bulk material is discretized with 3D volume elements while the fracture process zone is represented by 2D surface elements. In traditional finite-element frameworks, managing the connectivity between these topologically distinct domains and enforcing traction-separation relations requires dedicated mesh data structures and specialized element formulations.

In the proposed energy-centric framework, we treat the cohesive interface as an additional surface energy contribution to the global potential. The total energy of the system is expressed as the sum of the bulk energies of the two bodies, defined on domains $\Omega_a$ and $\Omega_b$, and the fracture energy defined on the cohesive interface $\mathrm{S}_\text{coh}$:
\begin{equation} 
\Psi(\boldsymbol{u}) = \int_{\Omega_{\mathrm{a}}}\psi_\varepsilon(\boldsymbol{u}_a)~\mathrm{d\Omega} + \int_{\mathrm{\Omega_b}}\psi_\varepsilon(\boldsymbol{u}_b)~\mathrm{d\Omega} + \int_{\mathrm{S}_{\text{coh}}} \psi_{\text{coh}}(\delta)~\mathrm{dS}~,
\end{equation}
where the interface potential $\psi_\text{coh}$ depends on effective opening $\delta$, defined as the norm of  the displacement jump vector $[\![\boldsymbol{u}]\!]=\boldsymbol{u}_b-\boldsymbol{u}_a$ across the interface. 

We employ an exponential cohesive law governed by the fracture energy $\Gamma$ and the critical strength $\sigma_c$, with potential and traction response:
\begin{equation} 
\psi(\delta) = \Gamma \left[ 1 - \left(1 + \frac{\delta}{\delta_c}\right) \exp\left(-\frac{\delta}{\delta_c}\right) \right], \quad T(\delta) = \frac{\Gamma}{\delta_c^2} \delta \exp\left(-\frac{\delta}{\delta_c}\right)~, 
\end{equation}
where the characteristic length is given by $\delta_c=(\Gamma \exp(-1))/\sigma_c$. 

To ensure thermodynamic consistency, the fracture process must be irreversible. We enforce this condition directly at the level of the potential by tracking the maximum historical opening $\delta_\text{max}$. During unloading, the material response follows a linear path to the origin rather than retracing the softening curve:
\begin{equation} 
\psi_{\text{irrev}}(\delta) = 
\begin{cases} \psi(\delta) & \text{if } \delta = \delta_\text{max} \quad \text{(Loading)} \\ \psi(\delta_\text{max}) - \frac{1}{2}(\delta_\text{max} - \delta) \left( T(\delta_\text{max}) + \frac{T(\delta_\text{max})}{\delta_\text{max}} \delta \right) & \text{if } \delta < \delta_\text{max} \quad \text{(Unloading)} \end{cases} 
\end{equation}

We verify this formulation by simulating Mode-I fracture in a double cantilever beam under quasi-static loading, as shown in \cref{fig:cohesive-fracture}. We use the bulk material and interface parameters reported in ~\cite{pundir_transonic_2024}. We compute the displacement jump at the interface nodes, evaluate it at the quadrature points, and integrate the resulting cohesive energy density along the interfaces. The following code illustrates the implementation of cohesive energy, including the enforcement of irreversibility:

\begin{lstlisting}
def exponential_cohesive_energy(jump: Array, delta_max: float) -> float:
    delta = compute_opening(jump) # compute effective opening
    
    def opening_branch(delta):
        loading = lambda delta : exponential_potential(delta, Gamma, delta_c)

        def unloading(delta):
            ... # compute psi_irrev as in Equation 12
            return psi_irrev
        return jax.lax.cond(delta > delta_max, loading, unloading, delta)

    penetration = lamda delta : 0.5 * penalty * delta**2
    return jax.lax.cond(delta > 1e-8, opening_branch, penetration_branch, delta)

def total_cohesive_energy(u_flat: Array, delta_max: Array) -> float:
    u = u_flat.reshape(-1, n_dofs_per_node)
    jump = u.at[top_interface_nodes, :].get() - u.at[bottom_interface_nodes, :].get()
    jump_quad = op_quad.eval(jump)   # Operator defined at the interface
    cohesive_energy_density = .... # call exponential_cohesive_energy
    return op_quad.integrate(cohesive_energy_density)
\end{lstlisting}

We solve the resulting nonlinear problem using a matrix-free Newton--Krylov method. As a validation, we compare the total dissipated energy with the theoretical value $\Gamma\times A_\text{crack}$ required to fully separate the interface. The close agreement observed in \Cref{fig:cohesive-fracture} confirms that the mixed-dimensional coupling correctly transmits forces and dissipates energy, without relying on specialized ``interface element'' classes.

\begin{figure}[!t]
    \centering
    \includegraphics[width=\linewidth]{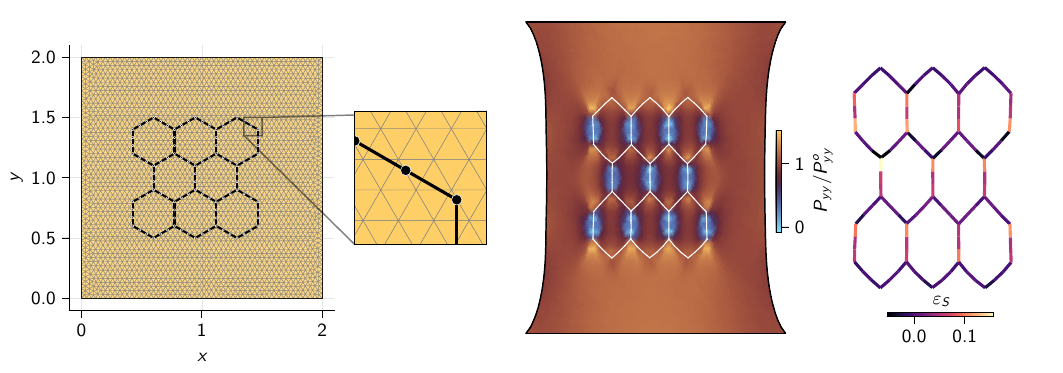}
    \caption{Fiber-reinforced soft material under uniaxial tension. (left) Schematic of a hexagonal fiber network embedded in a compressible neo-Hookean material. The two domains are discretized non-conformingly, such that fiber nodes are embedded within triangular bulk elements (see inset). The fibers have cross-sectional area $A_\text{fiber}$ and stiffness $E_\text{fiber} = 100 E_\text{bulk}$. (middle) Deformed bulk configuration at an applied stretch of $1.2$, colored by the first Piola–Kirchhoff stress component $P_{yy}$ normalized by $P_{yy}^{o}$, which denotes the bulk stress under the same imposed deformation without fibers. (right) The deformed configuration of the embedded fiber network, colored by the strain $\varepsilon_S$ in each fiber.}
    \label{fig:embedded-fiber}
\end{figure}

\subsubsection{Embedded fiber network in a soft material}\label{sec:embedded}

As a second example of mixed-dimensional coupling, we consider a fiber-reinforced composite in which stiff 1D fibers are embedded in a soft 2D or 3D bulk material, see \Cref{fig:embedded-fiber}. Such configurations arise in biological tissues (\eg{}, muscle fibers) and engineering composites. A central difficulty in these problems is that the fiber geometry typically does not align with the bulk mesh. The fibers may intersect bulk elements at arbitrary positions and orientations.

We address this using an embedded-element approach formulated entirely at the level of the total potential energy. We treat the fiber and bulk contributions as separate energy terms. The bulk energy $\Psi_\mathrm{bulk}$ is integrated over the triangular elements, while the fiber energy $\Psi_{\text{fiber}}$ is integrated along 1D line elements. The total energy is thus given as  
\begin{equation} 
\Psi(\boldsymbol{u}) = \underbrace{\int_{\Omega_{\text{bulk}}} \psi_\varepsilon(\nabla \boldsymbol{u}) ~\mathrm{d\Omega}}_{\Psi_\mathrm{bulk}} + \underbrace{\int_{\Gamma_{\text{fiber}}} \psi_{\text{fiber}}(\epsilon_{\text{fiber}})~ \mathrm{dS}}_{\Psi_{\text{fiber}}}~,
\end{equation}
where we define the fiber strain energy density as $\psi_\text{fiber}(\varepsilon_\text{fiber}) = \frac{1}{2}E_\text{fiber}A_\text{fiber} \varepsilon_\text{fiber}^2$. In this example, we model the bulk as a compressible neo-Hookean material.

The key step is coupling the fiber kinematics to the bulk deformation. Because the meshes are non-conforming, we perform a geometric search prior to the simulation. For each fiber node, we identify the containing triangular bulk element and compute its barycentric coordinates within that element. During the simulation, we use this precomputed mapping to interpolate the bulk displacement field to the fiber nodes. This procedure effectively ``ties'' the fiber to the surrounding bulk without introducing additional DoFs or modifying the mesh topology. 

We then compute the fiber strain at the quadrature-points and the strain energy of the fiber based on this interpolated displacement as:

\begin{lstlisting}
def fiber_strain_energy(u_flat: Array, barycenters: Array) -> float:
    u = u_flat.reshape(-1, n_dofs_per_node) # bulk displacments
    u_at_nodes = ... # project the bulk displacement to the fiber nodes using barycenters
    u_grad = op_line.grad(u_at_nodes) 
    energy_density = compute_fiber_energy(u_grad)
    return op_line.integrate(energy_density)
\end{lstlisting}

We finally add this scalar contribution to the global energy functional. AD propagates the fiber stiffness consistently to the bulk DoFs through the chain rule. As a result, the framework handles non-conforming embedded reinforcements in a purely variational manner, without specialized element types or intrusive modifications to the discretization.

\begin{figure}[!t]
	\centering
	\includegraphics[width=\linewidth]{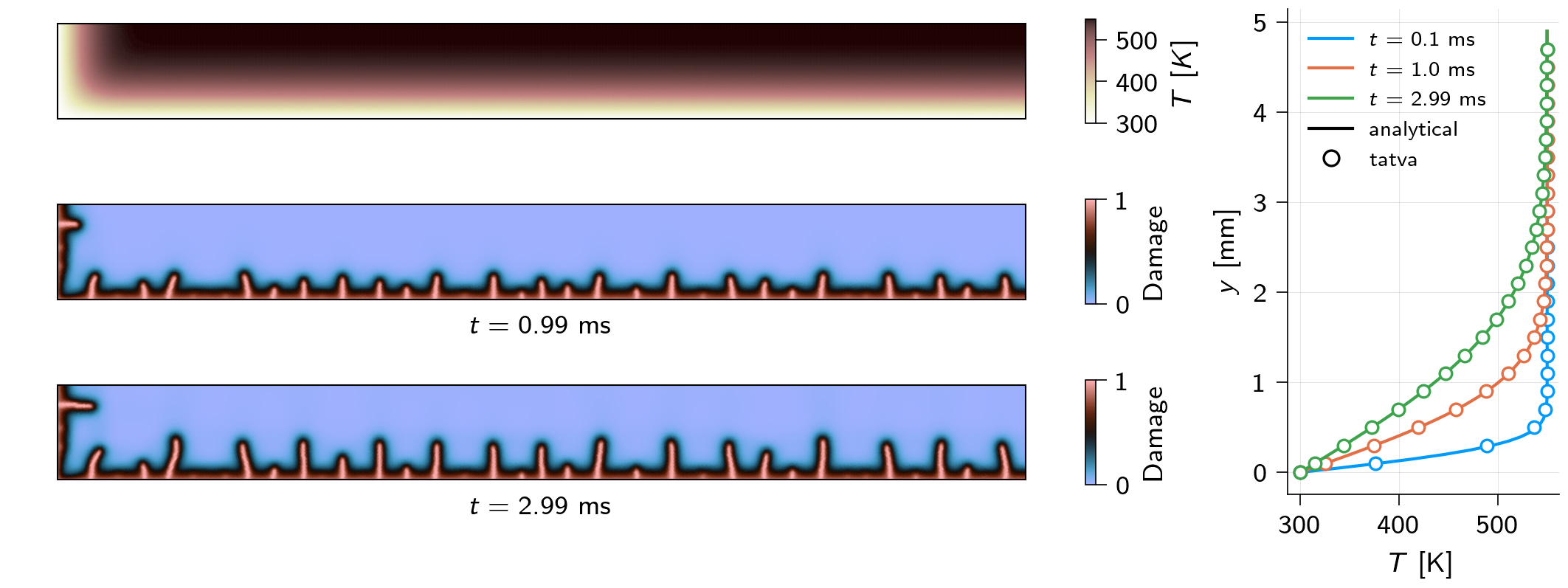}
	\caption{Phase-field fracture during thermal shock of a ceramic slab ($50\times4.9$ mm) initially at $T=550$ K, subjected to sudden cooling at the boundaries to $T=300$ K. The material and geometric parameters for the problem are adopted from \cite{mandal_fracture_2021}. (left-top) Temperature distribution at $t=2.99$ ms and (left-middle-bottom) corresponding crack patterns at $t=0.99$ ms and $2.99$ ms. Cracks initiate at the cooled surface and propagate inward under thermal gradient, consistent with experimental findings~\cite{shao_effect_2011}. (right) Temperature profile along the slab height at three time steps compared with the analytical solution~\cite{mandal_fracture_2021}, showing excellent agreement. }
	\label{fig:thermal_shock}
\end{figure}

\subsection{Multiphysical coupling}\label{sec:ex-thermal-shock}

Multiphysics problems involve the interaction of multiple field variables within a single model. In a variational setting, these interactions can be expressed through a unified energy functional that aggregates the contributions of all physical processes. We demonstrate the capability of the proposed framework to handle such couplings by simulating thermal shock fracture in a 2D ceramic slab~\cite{shao_effect_2011}.

This problem involves a three-way coupling: the thermal field induces mechanical strains, the elastic energy drives fracture, and the evolving phase-field degrades the material stiffness.

At each time step, we minimize an incremental potential $\Psi(\boldsymbol{u},d,T)$, which combines elastic-thermal energy, fracture surface energy, and a transient heat potential:
\begin{equation} 
\begin{split}
\Psi &= \underbrace{\int_{\Omega} \left[ g(d) \psi^{+}_\varepsilon(\boldsymbol{u}, \Delta T) + \psi^{-}_\varepsilon(\boldsymbol{u}, \Delta T) \right] ~ \mathrm{d\Omega}}_{\text{Elastic-Thermal Energy}} + \underbrace{\int_\Omega \Gamma \left[ \frac{(d-1)^2}{2l_0} + \frac{l_0}{2} |\nabla d|^2 \right] ~ \mathrm{d\Omega}}_{\text{Fracture Surface Energy}} \\ &+ \underbrace{\int_\Omega \frac{1}{2}\Delta t (T^{t}-T^{t-1})^2~\mathrm{d\Omega} + \int_\Omega \kappa \nabla T\cdot\nabla T~\mathrm{d\Omega}}_{\text{Transient Heat Potential}, \Psi_{\mathrm{th}}} ~.
\end{split}
\end{equation}
Here, $\psi^+_\varepsilon$ and $\psi^-_\varepsilon$ denote the tensile and compressive parts of the elastic strain energy density~\cite{Ambati2014}, including thermal expansion effects. Further, $g(d)=(1-d)^2$ is the degradation function, $\Gamma$ is the fracture energy, and $l_0$ is the phase-field length scale. 

The thermal contribution, $\Psi_\text{th}$, corresponds to a backward Euler discretization of the transient heat equation and is implemented as 

\begin{lstlisting}
def thermal_energy_density(grad_T, kappa):
    return kappa * jnp.vdot(grad_T, grad_T)

def total_thermal_energy(T_new: Array, T_prev: Array, dt: float) -> float:
    T_new_quad = op.eval(T_new)               # temperature at time t
    T_prev_quad = op.eval(T_prev)             # temeprature at time t-1
    grad_T = op.grad(T_new)
    diff_term = thermal_energy_density(grad_T, kappa)
    inertial_term = (0.5 / dt) * (T_new_quad - T_prev_quad)**2 # temporal
    return op.integrate(diff_term + inertial_term)
\end{lstlisting}  
The elastic and fracture energy terms follow directly from the formulations introduced in earlier examples.

For this simulation, we employed a staggered (operator-split) solution strategy. At each time step, we sequentially minimize the potential with respect to the displacement field $\boldsymbol{u}$, the phase-field $d$, and the temperature field $T$, iterating until convergence.

As shown in \Cref{fig:thermal_shock}, cracks initiate at the surface and propagate inward under the imposed thermal gradient, in qualitatively agreement with experimental observations~\cite{shao_effect_2011}. Furthermore, the temperature distribution inside the slab exhibits excellent agreement with the corresponding analytical solution~\cite{mandal_fracture_2021}.

A key advantage of the energy-centric architecture is that switching from a staggered to a fully monolithic solution strategy requires no modifications to the physics implementation. To perform a monolithic solve, we take a similar approach as in \cref{sec:lagrange} and concatenate all DoFs into a single global state vector $\boldsymbol{z}=[\boldsymbol{u},d,T]$ and define the total energy directly in terms of this vector:
\begin{lstlisting}
def total_energy(z_flat: Array) -> float:
    u, d, T = slice_and reshape_z(z_flat)
    ... # compute strain energy density
    ... # compute fracture energy density
    ... # compute thermal potential density
    return op.integrate(strain_energy_density) 
           + op.integrate(fracture_energy_density) 
           + op.integrate(thermal_energy_density)

monolithic_residual = jax.jacrev(total_energy)
monolithic_tangent = jax.jvp(total_energy, (z_prev,), (dz))[1]
\end{lstlisting}
AD then computes the full Hessian (or JVP) with respect to $\boldsymbol{z}$, automatically capturing all cross-coupling terms, such as $\nabla^2_{\boldsymbol{u},T}\Psi$, without manual block assembly.

\begin{figure}[!t]
	\centering
	\includegraphics[width=\linewidth]{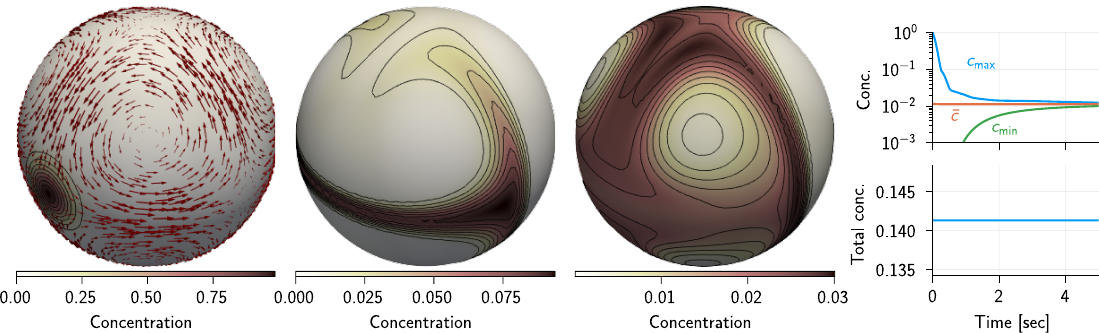}
	\caption{Advection-diffusion of a scalar concentration on a spherical surface. (left) Evolution of the concentration field at three time steps. The initial Gaussian pulse advects under a divergence-free velocity field (red arrows scaled by the velocity magnitude) and gradually diffuses. (middle) Concentration profile at $t=1$ s showing dominant advective transport followed by diffusion. (right-top) Evolution of maximum, mean, and minimum concentration values over time, illustrating homogenization due to diffusion. (right-bottom) Total surface concentration $\int c ~d\text{S}$ remains constant, confirming mass conservation and numerical robustness.}
	\label{fig:advection-diffusion}
\end{figure}

\subsection{Non-Variational problem on non-Euclidean manifolds}\label{sec:ex:non-variational}

To demonstrate that the framework is not limited to conservative systems, we solve an advection-diffusion problem for a scalar field $c$ on a spherical surface $\mathrm{S}$~\cite{rognes_automating_2013}. The governing equation reads
\begin{equation}
	\frac{\partial c}{\partial t} + \boldsymbol{u} \cdot \mathrm{\nabla_s} c = D \mathrm{\Delta_s} c~,
\end{equation}
where $D$ denotes the diffusivity, $\boldsymbol{u}$ is a velocity field tangent to the sphere, $\mathrm{\Delta_s}$ is the Laplace-Beltrami operator, and $\mathrm{\nabla_s}$ is the surface gradient.

Unlike previous examples, the advective term makes the operator non-symmetric. Consequently, no scalar potential $\Psi$ exists whose variation yields the governing equation. Instead, we formulate the problem using the principle of virtual work, as outlined in \cref{sec:non-variational}. The weak form is given as:
\begin{equation}
	\int_{\mathrm{S}} \frac{\partial c}{\partial t} v \, \mathrm{dS} + \int_{\mathrm{S}} (\boldsymbol{u} \cdot \mathrm{\nabla_s} c) v \, \mathrm{dS} + \int_{\mathrm{S}} D \mathrm{\nabla_s} c \cdot \mathrm{\nabla_s} v \, \mathrm{dS} = 0~,
\end{equation}
where $v$ denotes the test function. For this example, we restrict the velocity field to be \textit{divergence-free}, \ie{} $\mathrm{\nabla_s} \cdot \boldsymbol{u} = 0$, which simplified the weak formulation.

We implement this formulation using the same operators employed in the previous energy-based examples. No special plugins or problem-specific extensions are required. The scalar virtual-work functional is defined as

\begin{lstlisting}
def diffusion_energy_density(grad_c, grad_v, D):
    return D * jnp.vdot(grad_c, grad_v)

def total_virtual_work(c_flat, v_flat):
    grad_c = op.grad(c_flat)
    grad_v = op.grad(v_flat)
    density = diffusion_energy_density(grad_c, grad_v, D)
    return op.integrate(density)

compute_internal_force = jax.jacrev(total_virtual_work, argnums=1)
residual = compute_internal_force(u_flat,  v_flat=jnp.zeros(n_dofs))
\end{lstlisting}

AD then produces the residual vector and the action of the non-symmetric tangent operator within the same monolithic architecture.

We prescribe a velocity field corresponding to rigid-body rotation about the $z$-axis, $\boldsymbol{u} = \boldsymbol{\omega} \times \boldsymbol{x}$, which is analytically divergence-free on the sphere. Starting from a Gaussian concentration pulse (see \cref{fig:advection-diffusion}-left), the solution first advects along the imposed velocity field and subsequently diffuses over the surface, see \cref{fig:advection-diffusion}-middle. As shown in \cref{fig:advection-diffusion}-right, the maximum and minimum concentrations converge toward the mean concentration over time, while the total concentration remains constant, confirming mass conservation. 

This example highlights two key aspects of our framework. First, it supports non-variational and non-symmetric operators without modifying the architecture. Second, extending the formulation to a curved manifold requires only redefining the element-level Jacobian and gradient operators (see Appendix \ref{appendix:tri3manifold}). Once these geometric operators are provided, the same global energy/virtual-work structure and AD machinery apply unchanged.

\subsection{Integration of data-driven approaches} \label{sec:integrate_data_driven}

As a final demonstration of the proposed framework's versatility, we integrate data-driven models directly into the energy-centric finite-element formulation. Because the proposed architecture operates entirely on scalar functionals and relies on AD, it allows neural networks and other differentiable models to be incorporated without modifying the solver or assembly logic.

\subsubsection{Neural Constitutive Modeling (NCM)}\label{sec:ex:ncm}

NCMs offer a flexible way to represent complex material behavior beyond classical phenomenological laws. Recent works~\cite{latyshev_expressing_2025, alheit_commet_2026} have successfully integrated NCMs into finite-element simulations by coupling differentiable machine learning libraries with conventional solvers. However, these implementations are typically embedded within legacy, element-centric frameworks that rely on scatter-add assembly. As a result, while the neural network enhances constitutive expressivity, the overall architecture remains constrained by the underlying assembly workflow.

In the proposed framework, NCMs are not external components but native parts of the energy formulation. Because the entire framework is built on a differentiable infrastructure, a data-driven potential $\Psi_{NN}(\mathbf{F};\theta)$ is treated exactly like a classical strain energy density (see \cref{sec:ex1:bcs} and \cref{sec:ex2:rigid}). At the quadrature level, we simply replace the analytical energy density with a neural network evaluation,
\begin{lstlisting}
def strain_energy_density(grad_u):
    I = jnp.eye(3)
    F = I + grad_u 
    C = F.T @ F   # computing Cauchy strain
    ... # computing invariants
    return neural_network(I1, J) # calling neural network
\end{lstlisting}
The neural network maps strain invariants, such as $I_1 = \mathrm{tr}(C)$ and $J = \det(F)$, to a scalar strain energy density. Since the total potential energy is defined globally, AD propagates derivatives through the neural network and into the residual and tangent operators without any additional implementation effort.

\begin{figure}[!t]
	\centering
	\includegraphics[width=\linewidth]{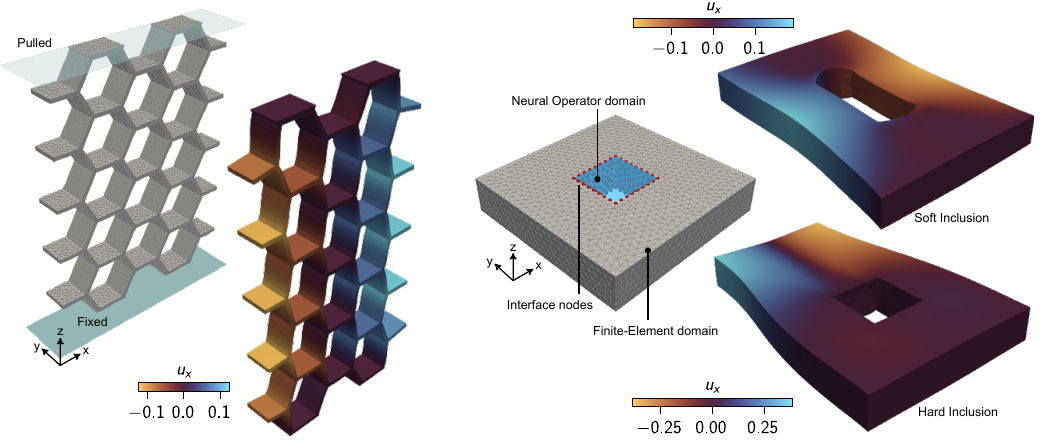}
	\caption{Incorporating neural networks within the energy-centric framework. (left) Hexagonal meta-material under uniaxial tension applied in the $z$-direction. A neural constitutive model is evaluated at each quadrature point to predict strain energy density from strain invariants. The color map shows the displacement component in the $x-$direction, illustrating a physically plausible deformation response. (right) Neural-Operator Element Method for multiscale modeling. The blue region indicates the subdomain represented by a neural operator, while the red nodes denote the interface between the neural domain and the classical FE domain. Loading is applied in the $y-$direction. The upper and lower subfigures correspond to soft and stiff neural inclusions, respectively. The color maps show the displacement component in the $x-$direction.}
	\label{fig:neural_constitutive_modeling}
\end{figure}

To illustrate this integration, we simulate a uniaxial tension test on a 3D hexagonal architected meta-material, see \cref{fig:neural_constitutive_modeling}-left. At each quadrature point, we evaluate a feed-forward multilayer perceptron with two hidden layers and \texttt{softplus} activation to ensure twice differentiability of the energy. To stabilize the initial iterations, we regularize the neural energy density as
\begin{equation} 
\psi_{\text{NN}}(I_1, J) = \left[ \text{NN}(I_1, J; \theta) - \text{NN}(3, 1; \theta) \right] + \psi_{\text{base}}(I_1, J)~, 
\end{equation}
where the shift enforces zero energy in the reference configuration and the convex base energy density ensures a positive-definite tangent stiffness.

For demonstration purposes, we employ an untrained synthetic network. The Newton solver converges in five iterations and produces physically plausible displacements (see \cref{fig:neural_constitutive_modeling}-left). This example highlights that neural constitutive models can be integrated seamlessly into the global energy formulation of our proposed framework, without modifying the solver, assembly logic, or differentiation pipeline. 

We emphasize that the framework is not limited to multilayer perceptrons. Any differentiable architecture supported by the AD backend, including Convolutional Neural Networks (CNNs), input convex neural networks (ICNNs)~\cite{thakolkaran_nn-euclid_2022}, constitutive artificial neural networks (CANNs)~\cite{thakolkaran_can_2025}, or input convex Kolmogorov-Arnold networks (ICKANs)~\cite{peirlinck_automated_2024}, can be incorporated in the same manner.

\subsubsection{Neural-Operator Element Method (NOEM)}\label{sec:ex:noem}

While neural constitutive modeling replaces local material behavior at the quadrature level, the same energy-centric architecture also supports domain-level learning approaches such as the NOEM~\cite{wang_accelerating_2025, wang_time_2025}. In NOEM, an entire subdomain, \eg{}, a region with intricate microstructure or geometric complexity, is replaced by a single neural element. Instead of discretizing this region explicitly, a pre-trained neural operator predicts its total potential energy from the displacement field prescribed on its boundary. 

A key challenge in conventional implementations is enforcing displacement continuity at the interface and deriving the stiffness terms that couple the neural region to the surrounding finite-element mesh. In the proposed framework, this coupling follows directly from energy superposition. We express the global potential as
\begin{equation} 
\Psi(\boldsymbol{u}) = \underbrace{\int_{\Omega_{\text{FEM}}} \psi_\varepsilon(\nabla \boldsymbol{u})~ \mathrm{d\Omega}}_{\Psi_{\mathrm{FEM}}} + \Psi_{\text{NO}}(\boldsymbol{u}_{\text{interface}}) - \Psi_{\text{ext}}~. 
\end{equation}
Here, $\Psi_\text{NO}$ is a scalar function that maps the interface displacements $\boldsymbol{u}_\text{interface}$ to the total energy of the neural inclusion. Because we define the neural contribution as an additive term in the global energy function, AD generates the interface forces and all cross-coupling stiffness terms without additional derivatives.

The implementation only requires extracting the interface displacements and passing them to the neural operator:

\begin{lstlisting}
def total_energy(u_flat: Array) -> float:
    u = u_flat.reshape(-1, n_dofs_per_node)
    energy_fem_domain = total_strain_energy(u_flat) # same function in Listing 1
    u_interface = u[interface_idx]   # extract displacements for interface nodes
    energy_neural_domain = neural_operator(u_interface) # psi_no

    return energy_fem_domain + energy_neural_domain
\end{lstlisting}

To demonstrate this formulation, we consider a 3D neural inclusion benchmark (see \cref{fig:neural_constitutive_modeling}-right). We generate an unstructured tetrahedral mesh of a cuboid domain containing a central square cutout. We do not mesh the cutout region. Instead, we treat it as a single neural element whose response depends only on the interface displacements. For this demonstration, we use an untrained multilayer perceptron with random weights to verify that the solver correctly handles domain coupling and derivative propagation.

To ensure well-posedness, we enforce zero energy at zero displacement by shifting the neural potential
\begin{equation} 
\Psi_{\text{NO}}(\boldsymbol{u}_{\text{interface}}) = \left[ \text{NN}(\boldsymbol{u}_{\text{interface}}; \theta) - \text{NN}(\mathbf{0}; \theta) \right] ~. 
\end{equation}
Because the neural inclusion contributes directly to the global energy, AD produces all interface forces and stiffness contributions consistently. The Newton-Raphson solver converges to a residual tolerance of $10^{-12}$ without manually deriving interface traction expressions or assembling dedicated coupling blocks.

Although we use an untrained synthetic network in this example, the same architecture applies to trained neural operators such as Fourier Neural Operator (FNO)~\cite{nguyen_universal_2025} or DeepONet~\cite{wang_accelerating_2025}. In practice, a trained operator can replace a geometrically complex domain while preserving a monolithic and fully differentiable formulation.

\section{Discussion} \label{sec:discussion}

In this work, we introduced an energy-centric finite-element framework that applies AD to a global potential functional to derive residual vectors and tangent operators. The central objective was to reconcile two traditionally conflicting properties in AD-based finite-element methods: high-level expressiveness at the level of the global variational formulation and computational scalability on modern hardware.

The examples presented in \cref{sec:examples} demonstrate that this reconciliation is possible within a single monolithic architecture. Across problems involving nonlinear materials, multi-point constraints, contact, mixed-dimensional couplings, multiphysics interactions, non-variational formulations, and data-driven models, the implementation followed the same structural pattern: define a scalar functional and differentiate it globally. No new assembly routines, element types, or coupling blocks were required. This illustrates that the versatility of the proposed framework does not arise from adding specialized modules, but from delegating derivative construction to AD at the level of the global energy.

By incorporating global AD, we shift the construction of residuals and tangent operators from manual assembly logic to the compiler. This enables two complementary strategies: JVP for matrix-free schemes and graph-coloring-based sparse differentiation for explicit matrix construction. Unlike conventional scatter-add assembly architectures, which remain element-centric, the proposed formulation operates directly on the global functional. As shown in \cref{sec:performance}, both strategies achieve linear scaling with respect to the number of DoFs and maintain stable throughput on GPU hardware. The absence of unstructured atomic scatter operations removes a major source of memory contention in traditional assembly workflows.

A notable consequence of representing the system through a single global functional is a numerical duality between robustness and scalability. Because both sparse differentiation and matrix-free Jacobian–vector products derive from the same global energy, the choice of linearization strategy becomes a numerical decision rather than an architectural constraint. For constrained or ill-conditioned systems, such as those involving Lagrange multipliers, contact, or incompressibility, forming an explicit sparse stiffness matrix enables the use of robust direct solvers. In contrast, for very large systems where storing even a sparse matrix becomes prohibitive, the same formulation can be evaluated in a matrix-free Newton–Krylov scheme using Jacobian–vector products. In both cases, the underlying physics definition remains unchanged. This flexibility contrasts with traditional element-centric architectures, where matrix-free and assembled formulations often require separate implementation pathways.

While the matrix-free approach based on JVP is entirely physics-agnostic, sparse matrix assembly through graph coloring requires prior knowledge of the sparsity pattern. For standard mechanics and saddle-point problems, this pattern is determined by mesh connectivity and is therefore known \emph{a priori}, as shown in Appendix \ref{appendix:sparsity}. However, for problems with evolving connectivity, such as non-local interactions, identifying the sparsity pattern in advance becomes more challenging. Future work will focus on detecting sparsity patterns directly from the global energy functional, for example, through structured perturbations or randomized probing techniques. In addition, more advanced coloring strategies, including star or acyclic coloring, may further reduce the number of required evaluations for sparse assembly~\cite{montoison_revisiting_2025}. 

Regarding solver integration, the present framework focuses on residual and tangent construction while remaining fully compatible with existing solver libraries. The output of AD is either a sparse matrix or a linear operator, which can be passed directly to external libraries such as PETSc~\cite{avenue_petsc_nodate} or algebraic multigrid solvers. As a result, the proposed architecture does not require rewriting established linear or nonlinear solver infrastructures.

We also observed practical memory limits when reconstructing large, sparse stiffness matrices on a single device. On a single GPU, sparse differentiation became constrained at approximately five million DoFs. This limitation originates from the growth of the compiler's internal representation of the computational graph rather than from the algorithmic complexity of the formulation itself. As the problem size increases, the graph eventually exceeds internal serialization limits. An approach to overcome this limitation can be multi-device execution, which our energy-centric structure provides a natural pathway for. Because the global functional is expressed as a sum of independent element contributions, the computation is inherently suitable for data-parallel distribution. By partitioning the mesh and sharding state vectors across multiple accelerators, each device only processes a local portion of the computational graph. This strategy reduces per-device memory pressure while preserving linear scaling. Investigating distributed execution through mechanisms such as \texttt{jax.sharding} remains an important direction for future work.

%
These results indicate that finite-element formulations can be expressed directly at the level of the global variational problem without incurring the dense-operator overhead of naive global differentiation. By combining sparse differentiation and matrix-free Jacobian-vector products within a single architecture, the framework removes the reliance on scatter-add element assembly while preserving strong scaling behavior on modern accelerators.

More broadly, the framework shifts the implementation effort from hand-written derivative and assembly kernels to the level of a global functional and its differentiation pipeline. For matrix-free solvers, this makes the formulation largely independent of an explicit sparsity structure. When an assembled matrix is required, sparse differentiation still relies on a mesh-defined sparsity pattern, but the derivative evaluation and operator construction are handled systematically by AD. This separation provides a practical route to implement complex coupled formulations while keeping the numerical backend flexible.

These observations suggest that monolithic, compiler-accelerated formulations provide a viable alternative to traditional assembly-based finite-element architectures for large-scale computational mechanics. Rather than optimizing legacy element-centric workflows for new hardware, the proposed approach redefines the representation of the global problem itself, aligning mathematical abstraction and computational execution within a single differentiable framework.

\section{Conclusion} \label{sec:conclusion}

In this work, we resolve a persistent limitation in AD–based finite-element methods: the trade-off between high-level expressiveness and computational scalability. Previous approaches forced a choice between globally differentiated formulations with clear variational structure but prohibitive $\mathcal{O}(N^2)$ complexity, and element-wise assembly schemes that scale efficiently but constrain how complex physics can be expressed.

Our proposed framework reconciles these extremes. By formulating problems strictly at the level of a global energy functional, we preserve the full expressive power of variational mechanics, making it straightforward to represent diverse problems, ranging from mixed-dimensional coupling and multiphysics, to non-conservative systems and data-driven models. At the same time, we avoid the computational penalties traditionally associated with global differentiation. Through graph-coloring-based sparse differentiation and matrix-free Jacobian-vector products, we achieve linear $\mathcal{O}(N)$ scaling while maintaining compatibility with modern accelerator hardware. Unlike traditional scatter-add assembly, our approach eliminates atomic contention and sustains constant throughput on GPUs and TPUs.

We provide this framework as an open-source library, \tatva~\cite{tatva_library_2026}. More broadly, this work demonstrates that global AD and large-scale performance are not mutually exclusive. By delegating assembly and differentiation to the compiler, we decouple physical modeling from low-level implementation details. This shift enables a more versatile, scalable, and future-proof foundation for computational mechanics.

\myacknowledgements{
    DSK and MP acknowledge support from the Swiss National Science Foundation under the SNSF starting grant (TMSGI2\_211655). All authors acknowledge support from the ETH Research Grant (24-1 ETH-020). The authors also thank Aryan Sinha for useful discussions on neural networks and Daniel Rayneau-Kirkhope for his feedback on the manuscript.
}
\insertacknowledgements

\section*{CRediT authorship contribution statement}
\textbf{Mohit Pundir}: Conceptualization,  Methodology, Software, Investigation, Formal Analysis, Data Curation, Visualisation, Writing -- original draft.\\
\textbf{Flavio Lorez}: Methodology, Software, Investigation, Formal Analysis, Data Curation, Visualisation, Writing -- original draft.\\
\textbf{David S. Kammer}: Conceptualisation, Supervision, Writing -- Review and Editing, Funding acquisition.
\section*{Declaration of competing interest}
The authors declare that they have no known competing financial interests or personal relationships that could have appeared to influence the work reported in this paper.

\section*{Code Availability}
All examples used in this paper are available in the tatva documentation~\cite{tatva_docs_2026}. The library \tatva{} is released as an open-source project~\cite{tatva_library_2026}.

\newpage

\appendix

\section{Coloring approaches for FEM} \label{appendix:coloring}
\begin{figure}[!h]
    \centering
    \includegraphics[width=0.5\linewidth]{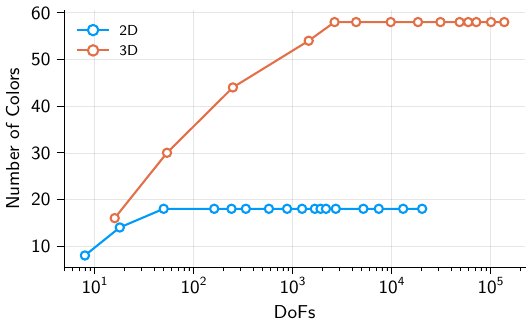}
    \caption{Number of colors as a function of total degrees of freedom for a 2D and a 3D mesh with 2 degrees of freedom per node. The element type for 2D is a 3-node triangle, and for 3D, it is a 4-node tetrahedral element.}
    \label{fig:coloring-max}
\end{figure}

In this work, we use the distance-2 greedy algorithm to color a tangent stiffness matrix. The coloring algorithm makes use of the sparsity pattern to identify and assign same color to columns whose non-zero entries do not share a row. The maximum number of colors required to color a matrix is governed by the mesh connectivity and the element type. For 2D element types with local support, such as 3-node triangular elements or 4-node quadrilateral elements, the maximum number of colors needed typically lies in the range of 30–40, regardless of the number of DoFs, as illustrated in \Cref{fig:coloring-max}. For 3D element types such as 4-node tetrahedral elements or 8-node hexahedron elements, typically about 80–90 colors are needed. This is because a single node in a 3D mesh is connected to more elements than in a 2D mesh. As illustrated in \Cref{fig:coloring-max}, these upper bounds are essentially independent of the overall problem size.

\section{Constructing sparsity pattern} \label{appendix:sparsity}
A sparsity pattern for a problem can be constructed in various ways. For single physical field problems, the standard approach is to use the mesh connectivity to construct the sparse pattern. Depending on the support (local or non-local) of the element type, one can detect which entries of the tangent stiffness matrix will be non-zero. For a Lagrangian functional $\mathcal{L}(\boldsymbol{u}, \lambda) = \Psi( \boldsymbol{u}) + \boldsymbol{\lambda} \cdot \boldsymbol{g}(\boldsymbol{u})$ such as in \Cref{sec:lagrange} where we enforce constraints $\boldsymbol{g}(\boldsymbol{u})$ via Lagrange multipliers, the tangent stiffness matrix has a characteristic block structure:
\begin{equation}
    \begin{bmatrix}
        \mathbf{K} & \mathbf{B}^\text{T} \\
        \mathbf{B} & \mathbf{0}
    \end{bmatrix}
\end{equation}
where $\mathbf{B}$  is given as
$$
\mathbf{B} = \frac{\partial \boldsymbol{g}}{\partial \mathbf{u}} ~.
$$

We see that the sparsity pattern of the above system is the union of the sparsity patterns of the matrices $\mathbf{K}$ and $\mathbf{B}$. The sparsity pattern of $\mathbf{K}$ can be constructed from mesh connectivity, while the sparsity pattern of $\mathbf{B}$ is computed using AD of the $g(\boldsymbol{u})$ function as 
\begin{lstlisting}
B = jax.jacrev(g)(jnp.zeros(n_dofs))
\end{lstlisting}

Thus, using these two sparsity patterns, one can construct the sparsity pattern for a Lagrangian functional.

\section{Implementation of triangular element in 3D space}\label{appendix:tri3manifold}

For the implementation of a 2D triangular element in 3D space, we simply redefine the Jacobian calculation and the gradient calculation to enable surface gradients in 3D space. 

\begin{lstlisting}
class Tri3Manifold(element.Tri3):
    """A 3-node linear triangular element on a 2D manifold embedded in 3D space."""

    def get_jacobian(self, xi: Array, nodal_coords: Array) -> tuple[Array, Array]:
        dNdr = self.shape_function_derivative(xi)
        J = dNdr @ nodal_coords  # shape (2, 2) or (2, 3)
        G = J @ J.T  # shape (2, 2)
        detJ = safe_sqrt(jnp.linalg.det(G))
        return J, detJ

    def gradient(self, xi: Array, nodal_values: Array, nodal_coords: Array) -> Array:
        dNdr = self.shape_function_derivative(xi)  # shape (2, 3)
        J, _ = self.get_jacobian(xi, nodal_coords)  # shape (2, 3)
        G_inv = jnp.linalg.inv(J @ J.T)  # shape (2, 2)
        J_plus = J.T @ G_inv  # shape (3, 2)
        dudxi = dNdr @ nodal_values  # shape (2, n_values)
        return J_plus @ dudxi  # shape (3, n_values)
\end{lstlisting}

\bibliographystyle{\mybibliostyle}

\begin{thebibliography}{55}
\providecommand{\natexlab}[1]{#1}
\providecommand{\url}[1]{\texttt{#1}}
\expandafter\ifx\csname urlstyle\endcsname\relax
  \providecommand{\doi}[1]{doi: #1}\else
  \providecommand{\doi}{doi: \begingroup \urlstyle{rm}\Url}\fi

\bibitem[Zienkiewicz and Taylor(2000)]{zienkiewicz_finite_2000}
O.~C. Zienkiewicz and Robert~L. Taylor.
\newblock \emph{The finite element method}.
\newblock Butterworth-Heinemann, Oxford ; Boston, 5th ed edition, 2000.
\newblock ISBN 978-0-7506-5049-6 978-0-7506-5055-7 978-0-7506-5050-2.

\bibitem[Bathe(1996)]{bathe_finite_1996}
Klaus-Jürgen Bathe.
\newblock \emph{Finite element procedures}.
\newblock Prentice Hall, Englewood Cliffs, New Jersey, second edition edition,
  1996.
\newblock ISBN 978-0-9790049-5-7.

\bibitem[Boffi et~al.(2013)Boffi, Brezzi, and Fortin]{boffi_mixed_2013}
Daniele Boffi, Franco Brezzi, and Michel Fortin.
\newblock \emph{Mixed {Finite} {Element} {Methods} and {Applications}},
  volume~44 of \emph{Springer {Series} in {Computational} {Mathematics}}.
\newblock Springer, Berlin, Heidelberg, 2013.
\newblock ISBN 978-3-642-36518-8 978-3-642-36519-5.
\newblock \doi{10.1007/978-3-642-36519-5}.
\newblock URL \url{https://link.springer.com/10.1007/978-3-642-36519-5}.

\bibitem[Margossian(2019)]{margossian_review_2019}
Charles~C. Margossian.
\newblock A review of automatic differentiation and its efficient
  implementation.
\newblock \emph{WIREs Data Mining and Knowledge Discovery}, 9\penalty0
  (4):\penalty0 e1305, 2019.
\newblock ISSN 1942-4795.
\newblock \doi{10.1002/widm.1305}.
\newblock URL \url{https://onlinelibrary.wiley.com/doi/abs/10.1002/widm.1305}.

\bibitem[Korelc and Wriggers(2016)]{korelc_automation_2016}
Jože Korelc and Peter Wriggers.
\newblock \emph{Automation of {Finite} {Element} {Methods}}.
\newblock Springer International Publishing, Cham, 2016.
\newblock ISBN 978-3-319-39003-1 978-3-319-39005-5.
\newblock \doi{10.1007/978-3-319-39005-5}.
\newblock URL \url{http://link.springer.com/10.1007/978-3-319-39005-5}.

\bibitem[Baydin et~al.(2017)Baydin, Pearlmutter, Radul, and
  Siskind]{baydin_automatic_2017}
Atılım~Günes Baydin, Barak~A. Pearlmutter, Alexey~Andreyevich Radul, and
  Jeffrey~Mark Siskind.
\newblock Automatic differentiation in machine learning: a survey.
\newblock \emph{J. Mach. Learn. Res.}, 18\penalty0 (1):\penalty0 5595--5637,
  January 2017.
\newblock ISSN 1532-4435.
\newblock URL \url{https://dl.acm.org/doi/10.5555/3122009.3242010}.

\bibitem[Rothe and Hartmann(2015)]{rothe_automatic_2015}
Steffen Rothe and Stefan Hartmann.
\newblock Automatic differentiation for stress and consistent tangent
  computation.
\newblock \emph{Archive of Applied Mechanics}, 85\penalty0 (8):\penalty0
  1103--1125, August 2015.
\newblock ISSN 1432-0681.
\newblock \doi{10.1007/s00419-014-0939-6}.
\newblock URL \url{https://doi.org/10.1007/s00419-014-0939-6}.

\bibitem[Schoenholz and Cubuk(2020)]{schoenholz_jax_2020}
Samuel~S. Schoenholz and Ekin~D. Cubuk.
\newblock {JAX}, {M}.{D}.: {A} {Framework} for {Differentiable} {Physics},
  December 2020.
\newblock URL \url{http://arxiv.org/abs/1912.04232}.
\newblock arXiv:1912.04232 [cond-mat, physics:physics, stat].

\bibitem[Carrer et~al.(2024)Carrer, Cezar, Bore, Ledum, and
  Cascella]{carrer_learning_2024}
Manuel Carrer, Henrique~Musseli Cezar, Sigbjørn~Løland Bore, Morten Ledum,
  and Michele Cascella.
\newblock Learning {Force} {Field} {Parameters} from {Differentiable}
  {Particle}-{Field} {Molecular} {Dynamics}.
\newblock \emph{Journal of Chemical Information and Modeling}, July 2024.
\newblock ISSN 1549-9596.
\newblock \doi{10.1021/acs.jcim.4c00564}.
\newblock URL \url{https://doi.org/10.1021/acs.jcim.4c00564}.

\bibitem[Toshev et~al.(2024)Toshev, Ramachandran, Erbesdobler, Galletti,
  Brandstetter, and Adams]{toshev_jax-sph_2024}
Artur~P. Toshev, Harish Ramachandran, Jonas~A. Erbesdobler, Gianluca Galletti,
  Johannes Brandstetter, and Nikolaus~A. Adams.
\newblock {JAX}-{SPH}: {A} {Differentiable} {Smoothed} {Particle}
  {Hydrodynamics} {Framework}, March 2024.
\newblock URL \url{http://arxiv.org/abs/2403.04750}.
\newblock arXiv:2403.04750 [physics].

\bibitem[Xue et~al.(2023)Xue, Liao, Gan, Park, Xie, Liu, and
  Cao]{xue_jax-fem_2023}
Tianju Xue, Shuheng Liao, Zhengtao Gan, Chanwook Park, Xiaoyu Xie, Wing~Kam
  Liu, and Jian Cao.
\newblock {JAX}-{FEM}: {A} differentiable {GPU}-accelerated {3D} finite element
  solver for automatic inverse design and mechanistic data science.
\newblock \emph{Computer Physics Communications}, 291:\penalty0 108802, October
  2023.
\newblock ISSN 0010-4655.
\newblock \doi{10.1016/j.cpc.2023.108802}.
\newblock URL
  \url{https://www.sciencedirect.com/science/article/pii/S0010465523001479}.

\bibitem[Bleyer(2024)]{bleyer_numerical_2024}
Jeremy Bleyer.
\newblock Numerical tours of {Computational} {Mechanics} with {FEniCSx},
  January 2024.
\newblock URL \url{https://doi.org/10.5281/zenodo.10470942}.

\bibitem[Latyshev et~al.(2025)Latyshev, Bleyer, Maurini, and
  Hale]{latyshev_expressing_2025}
Andrey Latyshev, Jérémy Bleyer, Corrado Maurini, and Jack Hale.
\newblock Expressing general constitutive models in {FEniCSx} using external
  operators and algorithmic automatic differentiation.
\newblock \emph{Journal of Theoretical, Computational and Applied Mechanics},
  September 2025.
\newblock ISSN 2726-6141.
\newblock \doi{10.46298/jtcam.14449}.
\newblock URL \url{https://jtcam.episciences.org/16616}.

\bibitem[Vigliotti and Auricchio(2021)]{vigliotti_automatic_2021}
Andrea Vigliotti and Ferdinando Auricchio.
\newblock Automatic {Differentiation} for {Solid} {Mechanics}.
\newblock \emph{Archives of Computational Methods in Engineering}, 28\penalty0
  (3):\penalty0 875--895, May 2021.
\newblock ISSN 1886-1784.
\newblock \doi{10.1007/s11831-019-09396-y}.
\newblock URL \url{https://doi.org/10.1007/s11831-019-09396-y}.

\bibitem[Wu(2023)]{wu_framework_2023}
Gaoyuan Wu.
\newblock A framework for structural shape optimization based on automatic
  differentiation, the adjoint method and accelerated linear algebra.
\newblock \emph{Structural and Multidisciplinary Optimization}, 66\penalty0
  (7):\penalty0 151, July 2023.
\newblock ISSN 1615-147X, 1615-1488.
\newblock \doi{10.1007/s00158-023-03601-0}.
\newblock URL \url{http://arxiv.org/abs/2211.15409}.
\newblock arXiv:2211.15409 [cs, math].

\bibitem[Ichihara(2026)]{ichihara_naruki-ichiharafeax_2026}
Naruki Ichihara.
\newblock Naruki-{Ichihara}/feax, January 2026.
\newblock URL \url{https://github.com/Naruki-Ichihara/feax}.
\newblock original-date: 2025-07-09T05:45:18Z.

\bibitem[Logg et~al.(2012)Logg, Ølgaard, Rognes, and Wells]{logg_ffc_2012}
Anders Logg, Kristian~B. Ølgaard, Marie~E. Rognes, and Garth~N. Wells.
\newblock {FFC}: the {FEniCS} {Form} {Compiler}.
\newblock In Anders Logg, Kent-Andre Mardal, and Garth~N. Wells, editors,
  \emph{Automated {Solution} of {Differential} {Equations} by the {Finite}
  {Element} {Method}}, volume~84 of \emph{Lecture {Notes} in {Computational}
  {Science} and {Engineering}}. Springer, 2012.
\newblock Section: 11.

\bibitem[Baratta et~al.(2023)Baratta, Dean, Dokken, Habera, Hale, Richardson,
  Rognes, Scroggs, Sime, and Wells]{baratta_dolfinx_2023}
Igor~A. Baratta, Joseph~P. Dean, Jørgen~S. Dokken, Michal Habera, Jack~S.
  Hale, Chris~N. Richardson, Marie~E. Rognes, Matthew~W. Scroggs, Nathan Sime,
  and Garth~N. Wells.
\newblock {DOLFINx}: the next generation {FEniCS} problem solving environment,
  2023.
\newblock Published: preprint.

\bibitem[Anderson et~al.(2021)Anderson, Andrej, Barker, Bramwell, Camier,
  Cerveny, Dobrev, Dudouit, Fisher, Kolev, Pazner, Stowell, Tomov, Akkerman,
  Dahm, Medina, and Zampini]{anderson_mfem_2021}
Robert Anderson, Julian Andrej, Andrew Barker, Jamie Bramwell, Jean-Sylvain
  Camier, Jakub Cerveny, Veselin Dobrev, Yohann Dudouit, Aaron Fisher, Tzanio
  Kolev, Will Pazner, Mark Stowell, Vladimir Tomov, Ido Akkerman, Johann Dahm,
  David Medina, and Stefano Zampini.
\newblock {MFEM}: {A} modular finite element methods library.
\newblock \emph{Computers \& Mathematics with Applications}, 81:\penalty0
  42--74, January 2021.
\newblock ISSN 0898-1221.
\newblock \doi{10.1016/j.camwa.2020.06.009}.
\newblock URL
  \url{https://www.sciencedirect.com/science/article/pii/S0898122120302583}.

\bibitem[Kaczmarczyk et~al.(2020)Kaczmarczyk, Ullah, Lewandowski, Meng, Zhou,
  Athanasiadis, Nguyen, Chalons-Mouriesse, Richardson, Miur, Shvarts, Wakeni,
  and Pearce]{kaczmarczyk_mofem_2020}
Łukasz Kaczmarczyk, Zahur Ullah, Karol Lewandowski, Xuan Meng, Xiao-Yi Zhou,
  Ignatios Athanasiadis, Hoang Nguyen, Christophe-Alexandre Chalons-Mouriesse,
  Euan~J. Richardson, Euan Miur, Andrei~G. Shvarts, Mebratu Wakeni, and
  Chris~J. Pearce.
\newblock {MoFEM}: {An} open source, parallel finite element library.
\newblock \emph{Journal of Open Source Software}, 5\penalty0 (45):\penalty0
  1441, January 2020.
\newblock ISSN 2475-9066.
\newblock \doi{10.21105/joss.01441}.
\newblock URL \url{https://joss.theoj.org/papers/10.21105/joss.01441}.

\bibitem[Richart et~al.(2024)Richart, Anciaux, Gallyamov, Frérot, Kammer,
  Pundir, Vocialta, Ramos, Corrado, Müller, Barras, Zhang, Ferry, Durussel,
  and Molinari]{richart_akantu_2024}
Nicolas Richart, Guillaume Anciaux, Emil Gallyamov, Lucas Frérot, David
  Kammer, Mohit Pundir, Marco Vocialta, Aurelia~Cuba Ramos, Mauro Corrado,
  Philip Müller, Fabian Barras, Shenghan Zhang, Roxane Ferry, Shad Durussel,
  and Jean-François Molinari.
\newblock Akantu: an {HPC} finite-element library for contact and dynamic
  fracture simulations.
\newblock \emph{Journal of Open Source Software}, 9\penalty0 (94):\penalty0
  5253, February 2024.
\newblock ISSN 2475-9066.
\newblock \doi{10.21105/joss.05253}.
\newblock URL \url{https://joss.theoj.org/papers/10.21105/joss.05253}.

\bibitem[Andrej et~al.(2025)Andrej, Kolev, and Lazarov]{andrej_scalable_2025}
Julian Andrej, Tzanio Kolev, and Boyan Lazarov.
\newblock Scalable {Analysis} and {Design} {Using} {Automatic}
  {Differentiation}, May 2025.
\newblock URL \url{http://arxiv.org/abs/2506.00746}.
\newblock arXiv:2506.00746 [math].

\bibitem[Pundir et~al.(2026)Pundir, Lorez, and Kammer]{tatva_library_2026}
M.~Pundir, F.~Lorez, and D.~S. Kammer.
\newblock tatva: An energy-centric automatic differentiation framework for
  finite-element analysis.
\newblock \url{https://github.com/smec-ethz/tatva}, 2026.
\newblock Open-source software.

\bibitem[Pundir and Kammer(2025)]{pundir_simplifying_2025}
Mohit Pundir and David~S. Kammer.
\newblock Simplifying {FFT}-based methods for solid mechanics with automatic
  differentiation.
\newblock \emph{Computer Methods in Applied Mechanics and Engineering},
  435:\penalty0 117572, February 2025.
\newblock ISSN 0045-7825.
\newblock \doi{10.1016/j.cma.2024.117572}.
\newblock URL
  \url{https://www.sciencedirect.com/science/article/pii/S0045782524008260}.

\bibitem[Montoison et~al.(2025)Montoison, Dalle, and
  Gebremedhin]{montoison_revisiting_2025}
Alexis Montoison, Guillaume Dalle, and Assefaw Gebremedhin.
\newblock Revisiting {Sparse} {Matrix} {Coloring} and {Bicoloring}, May 2025.
\newblock URL \url{http://arxiv.org/abs/2505.07308}.
\newblock arXiv:2505.07308 [math].

\bibitem[Gebremedhin et~al.(2005)Gebremedhin, Manne, and
  Pothen]{gebremedhin_what_2005}
Assefaw~Hadish Gebremedhin, Fredrik Manne, and Alex Pothen.
\newblock What {Color} {Is} {Your} {Jacobian}? {Graph} {Coloring} for
  {Computing} {Derivatives}.
\newblock \emph{SIAM Review}, 47\penalty0 (4):\penalty0 629--705, January 2005.
\newblock ISSN 0036-1445.
\newblock \doi{10.1137/S0036144504444711}.
\newblock URL \url{https://epubs.siam.org/doi/10.1137/S0036144504444711}.
\newblock Publisher: Society for Industrial and Applied Mathematics.

\bibitem[Curtis et~al.(1974)Curtis, Powell, and Reid]{curtis_estimation_1974}
A.~R. Curtis, M.~J.~D. Powell, and J.~K. Reid.
\newblock On the {Estimation} of {Sparse} {Jacobian} {Matrices}.
\newblock \emph{IMA Journal of Applied Mathematics}, 13\penalty0 (1):\penalty0
  117--119, February 1974.
\newblock ISSN 0272-4960.
\newblock \doi{10.1093/imamat/13.1.117}.
\newblock URL \url{https://doi.org/10.1093/imamat/13.1.117}.

\bibitem[Alheit et~al.(2026)Alheit, Peirlinck, and Kumar]{alheit_commet_2026}
Benjamin Alheit, Mathias Peirlinck, and Siddhant Kumar.
\newblock {COMMET}: {Orders}-of-magnitude speed-up in finite element method via
  batch-vectorized neural constitutive updates.
\newblock \emph{Computer Methods in Applied Mechanics and Engineering},
  452:\penalty0 118728, April 2026.
\newblock ISSN 0045-7825.
\newblock \doi{10.1016/j.cma.2026.118728}.
\newblock URL
  \url{https://www.sciencedirect.com/science/article/pii/S0045782526000022}.

\bibitem[Thakolkaran et~al.(2022)Thakolkaran, Joshi, Zheng, Flaschel,
  De~Lorenzis, and Kumar]{thakolkaran_nn-euclid_2022}
Prakash Thakolkaran, Akshay Joshi, Yiwen Zheng, Moritz Flaschel, Laura
  De~Lorenzis, and Siddhant Kumar.
\newblock {NN}-{EUCLID}: {Deep}-learning hyperelasticity without stress data.
\newblock \emph{Journal of the Mechanics and Physics of Solids}, 169:\penalty0
  105076, December 2022.
\newblock ISSN 0022-5096.
\newblock \doi{10.1016/j.jmps.2022.105076}.
\newblock URL
  \url{https://www.sciencedirect.com/science/article/pii/S0022509622002538}.

\bibitem[Joshi et~al.(2022)Joshi, Thakolkaran, Zheng, Escande, Flaschel,
  De~Lorenzis, and Kumar]{joshi_bayesian-euclid_2022}
Akshay Joshi, Prakash Thakolkaran, Yiwen Zheng, Maxime Escande, Moritz
  Flaschel, Laura De~Lorenzis, and Siddhant Kumar.
\newblock Bayesian-{EUCLID}: {Discovering} hyperelastic material laws with
  uncertainties.
\newblock \emph{Computer Methods in Applied Mechanics and Engineering},
  398:\penalty0 115225, August 2022.
\newblock ISSN 0045-7825.
\newblock \doi{10.1016/j.cma.2022.115225}.
\newblock URL
  \url{https://www.sciencedirect.com/science/article/pii/S0045782522003681}.

\bibitem[Han et~al.(2025)Han, Pundir, Fink, and Kammer]{han_learning_2025}
Zhichao Han, Mohit Pundir, Olga Fink, and David~S. Kammer.
\newblock Learning physics-consistent material behavior from dynamic
  displacements.
\newblock \emph{Computer Methods in Applied Mechanics and Engineering},
  443:\penalty0 118040, August 2025.
\newblock ISSN 0045-7825.
\newblock \doi{10.1016/j.cma.2025.118040}.
\newblock URL
  \url{https://www.sciencedirect.com/science/article/pii/S0045782525003123}.

\bibitem[Wang et~al.(2025{\natexlab{a}})Wang, Hakimzadeh, Ruan, and
  Goswami]{wang_accelerating_2025}
Wei Wang, Maryam Hakimzadeh, Haihui Ruan, and Somdatta Goswami.
\newblock Accelerating {Multiscale} {Modeling} with {Hybrid} {Solvers}:
  {Coupling} {FEM} and {Neural} {Operators} with {Domain} {Decomposition},
  April 2025{\natexlab{a}}.
\newblock URL \url{http://arxiv.org/abs/2504.11383}.
\newblock arXiv:2504.11383 [cs] version: 2.

\bibitem[Bradbury et~al.(2018)Bradbury, Frostig, Hawkins, Johnson, Leary,
  Maclaurin, Necula, Paszke, VanderPlas, Wanderman-Milne, and
  Zhang]{bradbury_jax_2018}
James Bradbury, Roy Frostig, Peter Hawkins, Matthew~James Johnson, Chris Leary,
  Dougal Maclaurin, George Necula, Adam Paszke, Jake VanderPlas, Skye
  Wanderman-Milne, and Qiao Zhang.
\newblock {JAX}: composable transformations of {Python}+{NumPy} programs, 2018.
\newblock URL \url{http://github.com/google/jax}.

\bibitem[Kidger and Garcia(2021)]{kidger_equinox_2021}
Patrick Kidger and Cristian Garcia.
\newblock Equinox: neural networks in {JAX} via callable {PyTrees} and filtered
  transformations, October 2021.
\newblock URL \url{http://arxiv.org/abs/2111.00254}.
\newblock arXiv:2111.00254 [cs].

\bibitem[Gerdes(2025)]{gerdes_mathisgerdesjax-autovmap_2025}
Mathis Gerdes.
\newblock mathisgerdes/jax-autovmap, September 2025.
\newblock URL \url{https://github.com/mathisgerdes/jax-autovmap}.
\newblock original-date: 2022-05-25T10:56:08Z.

\bibitem[noa()]{noauthor_openxlaxla_nodate}
openxla/xla: {A} machine learning compiler for {GPUs}, {CPUs}, and {ML}
  accelerators.
\newblock URL \url{https://github.com/openxla/xla}.

\bibitem[Pundir and Lorez(2026)]{tatva_docs_2026}
M.~Pundir and F.~Lorez.
\newblock tatva documentation and reproducible examples.
\newblock \url{https://smec-ethz.github.io/tatva-docs/}, 2026.
\newblock Documentation and executable examples (Google Colab).

\bibitem[BELYTSCHKO(2008)]{belytschko_finite_2008}
TED BELYTSCHKO.
\newblock The {Finite} {Element} {Method}: {Linear} {Static} and {Dynamic}
  {Finite} {Element} {Analysis}: {Thomas} {J}. {R}. {Hughes}.
\newblock \emph{Computer-Aided Civil and Infrastructure Engineering},
  4\penalty0 (3):\penalty0 245--246, November 2008.
\newblock ISSN 10939687.
\newblock \doi{10.1111/j.1467-8667.1989.tb00025.x}.
\newblock URL \url{http://doi.wiley.com/10.1111/j.1467-8667.1989.tb00025.x}.

\bibitem[Kirsch(1898)]{kirsch1898theorie}
Ernst~Gustav Kirsch.
\newblock Die theorie der elastizit t und die bed rfnisse der festigkeitslehre.
\newblock \emph{Zeitshrift des Vereines deutscher Ingenieure}, 42:\penalty0
  797--807, 1898.

\bibitem[Bleyer(2018)]{bleyer2018numericaltours}
Jeremy Bleyer.
\newblock \emph{Numerical Tours of Computational Mechanics with {FE}ni{CS}},
  2018.

\bibitem[Johnson(1987)]{johnson1987contact}
Kenneth~Langstreth Johnson.
\newblock \emph{Contact mechanics}.
\newblock Cambridge university press, 1987.

\bibitem[Popp(2012)]{Popp2012}
Alexander Popp.
\newblock Mortar {Methods} for {Computational} {Contact} {Mechanics} and
  {General} {Interface} {Problems}.
\newblock \penalty0 (September):\penalty0 34, 2012.
\newblock URL \url{http://mediatum.ub.tum.de/?id=1109994}.

\bibitem[Daversin-Catty et~al.(2021)Daversin-Catty, Richardson, Ellingsrud, and
  Rognes]{daversin-catty_abstractions_2021}
Cécile Daversin-Catty, Chris~N. Richardson, Ada~J. Ellingsrud, and Marie~E.
  Rognes.
\newblock Abstractions and {Automated} {Algorithms} for {Mixed} {Domain}
  {Finite} {Element} {Methods}.
\newblock \emph{ACM Trans. Math. Softw.}, 47\penalty0 (4):\penalty0
  31:1--31:36, September 2021.
\newblock ISSN 0098-3500.
\newblock \doi{10.1145/3471138}.
\newblock URL \url{https://dl.acm.org/doi/10.1145/3471138}.

\bibitem[Ortiz and Pandolfi(1999)]{Ortiz1999}
M.~Ortiz and A.~Pandolfi.
\newblock Finite-deformation irreversible cohesive elements for
  three-dimensional crack-propagation analysis.
\newblock \emph{International Journal for Numerical Methods in Engineering},
  44\penalty0 (9):\penalty0 1267--1282, 1999.
\newblock ISSN 0029-5981.
\newblock
  \doi{10.1002/(SICI)1097-0207(19990330)44:9<1267::AID-NME486>3.0.CO;2-7}.
\newblock URL
  \url{http://doi.wiley.com/10.1002/\%28SICI\%291097-0207\%2819990330\%2944\%3A9\%3C1267\%3A\%3AAID-NME486\%3E3.0.CO\%3B2-7}.
\newblock ISBN: 1097-0207.

\bibitem[Vocialta et~al.(2017)Vocialta, Richart, and Molinari]{Vocialta2017a}
M.~Vocialta, N.~Richart, and J.-F.~F. Molinari.
\newblock {3D} dynamic fragmentation with parallel dynamic insertion of
  cohesive elements.
\newblock \emph{International Journal for Numerical Methods in Engineering},
  109\penalty0 (12):\penalty0 1655--1678, March 2017.
\newblock ISSN 10970207.
\newblock \doi{10.1002/nme.5339}.
\newblock URL \url{http://doi.wiley.com/10.1002/nme.5339}.

\bibitem[Pundir et~al.(2024)Pundir, Adda-Bedia, and
  Kammer]{pundir_transonic_2024}
Mohit Pundir, Mokhtar Adda-Bedia, and David~S. Kammer.
\newblock Transonic and {Supershear} {Crack} {Propagation} {Driven} by
  {Geometric} {Nonlinearities}.
\newblock \emph{Physical Review Letters}, 132\penalty0 (22):\penalty0 226102,
  May 2024.
\newblock ISSN 0031-9007, 1079-7114.
\newblock \doi{10.1103/PhysRevLett.132.226102}.
\newblock URL \url{https://link.aps.org/doi/10.1103/PhysRevLett.132.226102}.

\bibitem[Mandal et~al.(2021)Mandal, Nguyen, Wu, Nguyen-Thanh, and
  de~Vaucorbeil]{mandal_fracture_2021}
Tushar~Kanti Mandal, Vinh~Phu Nguyen, Jian-Ying Wu, Chi Nguyen-Thanh, and Alban
  de~Vaucorbeil.
\newblock Fracture of thermo-elastic solids: {Phase}-field modeling and new
  results with an efficient monolithic solver.
\newblock \emph{Computer Methods in Applied Mechanics and Engineering},
  376:\penalty0 113648, April 2021.
\newblock ISSN 0045-7825.
\newblock \doi{10.1016/j.cma.2020.113648}.
\newblock URL
  \url{https://www.sciencedirect.com/science/article/pii/S0045782520308331}.

\bibitem[Shao et~al.(2011)Shao, Zhang, Xu, Zhou, Li, and Liu]{shao_effect_2011}
Yingfeng Shao, Yue Zhang, Xianghong Xu, Zhiliang Zhou, Wei Li, and Boyang Liu.
\newblock Effect of {Crack} {Pattern} on the {Residual} {Strength} of
  {Ceramics} {After} {Quenching}.
\newblock \emph{Journal of the American Ceramic Society}, 94\penalty0
  (9):\penalty0 2804--2807, 2011.
\newblock ISSN 1551-2916.
\newblock \doi{10.1111/j.1551-2916.2011.04728.x}.
\newblock URL
  \url{https://onlinelibrary.wiley.com/doi/abs/10.1111/j.1551-2916.2011.04728.x}.
\newblock \_eprint:
  https://ceramics.onlinelibrary.wiley.com/doi/pdf/10.1111/j.1551-2916.2011.04728.x.

\bibitem[Ambati et~al.(2014)Ambati, Gerasimov, and De~Lorenzis]{Ambati2014}
Marreddy Ambati, Tymofiy Gerasimov, and Laura De~Lorenzis.
\newblock A review on phase-field models of brittle fracture and a new fast
  hybrid formulation.
\newblock \emph{Computational Mechanics}, 55\penalty0 (2):\penalty0 383--405,
  2014.
\newblock ISSN 01787675.
\newblock \doi{10.1007/s00466-014-1109-y}.
\newblock ISBN: 0178-7675.

\bibitem[Rognes et~al.(2013)Rognes, Ham, Cotter, and
  McRae]{rognes_automating_2013}
M.~E. Rognes, D.~A. Ham, C.~J. Cotter, and A.~T.~T. McRae.
\newblock Automating the solution of {PDEs} on the sphere and other manifolds
  in {FEniCS} 1.2.
\newblock \emph{Geoscientific Model Development}, 6\penalty0 (6):\penalty0
  2099--2119, December 2013.
\newblock ISSN 1991-959X.
\newblock \doi{10.5194/gmd-6-2099-2013}.
\newblock URL \url{https://gmd.copernicus.org/articles/6/2099/2013/}.

\bibitem[Thakolkaran et~al.(2025)Thakolkaran, Guo, Saini, Peirlinck, Alheit,
  and Kumar]{thakolkaran_can_2025}
Prakash Thakolkaran, Yaqi Guo, Shivam Saini, Mathias Peirlinck, Benjamin
  Alheit, and Siddhant Kumar.
\newblock Can {KAN} {CANs}? {Input}-convex {Kolmogorov}-{Arnold} {Networks}
  ({KANs}) as hyperelastic constitutive artificial neural networks ({CANs}).
\newblock \emph{Computer Methods in Applied Mechanics and Engineering},
  443:\penalty0 118089, August 2025.
\newblock ISSN 0045-7825.
\newblock \doi{10.1016/j.cma.2025.118089}.
\newblock URL
  \url{https://www.sciencedirect.com/science/article/pii/S0045782525003615}.

\bibitem[Peirlinck et~al.(2024)Peirlinck, Linka, Hurtado, and
  Kuhl]{peirlinck_automated_2024}
Mathias Peirlinck, Kevin Linka, Juan~A. Hurtado, and Ellen Kuhl.
\newblock On automated model discovery and a universal material subroutine for
  hyperelastic materials.
\newblock \emph{Computer Methods in Applied Mechanics and Engineering},
  418:\penalty0 116534, January 2024.
\newblock ISSN 0045-7825.
\newblock \doi{10.1016/j.cma.2023.116534}.
\newblock URL
  \url{https://www.sciencedirect.com/science/article/pii/S0045782523006588}.

\bibitem[Wang et~al.(2025{\natexlab{b}})Wang, Hakimzadeh, Ruan, and
  Goswami]{wang_time_2025}
Wei Wang, Maryam Hakimzadeh, Haihui Ruan, and Somdatta Goswami.
\newblock Time {Marching} {Neural} {Operator} {FE} {Coupling}: {AI}
  {Accelerated} {Physics} {Modeling}, August 2025{\natexlab{b}}.
\newblock URL \url{http://arxiv.org/abs/2504.11383}.
\newblock arXiv:2504.11383 [cs].

\bibitem[Nguyen and Schneider(2025)]{nguyen_universal_2025}
Binh~Huy Nguyen and Matti Schneider.
\newblock Universal {Fourier} {Neural} {Operators} for {Micromechanics}, July
  2025.
\newblock URL \url{http://arxiv.org/abs/2507.12233}.
\newblock arXiv:2507.12233 [cs].

\bibitem[Avenue()]{avenue_petsc_nodate}
South~Cass Avenue.
\newblock {PETSc} {Users} {Manual}.

\end{thebibliography}
%

 %

%

\end{document}